\documentclass[11pt]{article}
\usepackage{amsmath,amsxtra}
\usepackage{bm}
\usepackage{amssymb}
\usepackage{mathrsfs}
\usepackage{amsfonts}
\usepackage{amscd}
\usepackage[dvipdf]{graphics}
\usepackage[dvips]{graphicx}
\usepackage{graphicx,subfigure}
\usepackage{latexsym}
\usepackage{epstopdf}
\usepackage{latexsym}
\usepackage{amsmath}
\newcommand{\beqa}{\begin{eqnarray}}
\newcommand{\eeqa}{\end{eqnarray}}
\newcommand{\ba}{\begin{eqnarray*}}
\newcommand{\ea}{\end{eqnarray*}}

\date{}

\newtheorem{definition}{Definition}

\def\div{{\,\rm div \,}}

\def\div{{\,\rm div \,}}

\def\be{{\begin{equation}}}
\def\ee{{\end{equation}}}
\def\Om{{\Omega}}
\def\na{{\nabla}}
\def\Ga{{\Gamma}}
\def\pl{{\partial}}
\def\epl{\epsilon}

\def\beq{\arraycolsep=1.5pt\begin{eqnarray}}
\def\eeq{\end{eqnarray}}

\newfont{\Blackboard}{msbm10 scaled 1200}

\newfont{\roma}{cmr10 scaled 1200}

\def\<{{\langle}}
\def\>{{\rangle}}

\newtheorem{thm}{{}\hskip\parindent Theorem}[section]
\newtheorem{lem}{{}\hskip\parindent Lemma}[section]

\newtheorem{cor}{{}\hskip\parindent Corollary}[section]

\newtheorem{rem}{{}\hskip\parindent Remark}[section]

\def\be{\begin{equation}}
\def\ee{\end{equation}}
\def\beq{\arraycolsep=1.5pt\begin{eqnarray}}
\def\eeq{\end{eqnarray}}

\large
\title{\bf Energy decay and global smooth solutions for a free boundary fluid-nonlinear elastic structure interface model with  boundary dissipation}
\date{}

\author{Yizhao Qin\quad
Peng-Fei Yao\thanks{Corresponding author.\ Email: pfyao@iss.ac.cn}\quad \\[0.3cm]
Key Laboratory of  Systems and Control\\
Institute of Systems Science,
Academy of Mathematics and Systems Science\\
Chinese Academy of Sciences, Beijing 100190, P. R. China\\
School of Mathematical Sciences\\
University of Chinese Academy of Sciences, Beijing 100049, China
}

\begin{document}
\maketitle
\footnote{This work is  supported by the National
Science Foundation of China, grants  no. 61473126 and no. 61573342, and Key Research Program of Frontier Sciences, CAS, no. QYZDJ-SSW-SYS011.}
\begin{quote}
\begin{small}
{\bf Abstract} \,\,\,We consider a nonlinear, free boundary fluid-structure interaction model in a bounded domain. The viscous incompressible fluid interacts with a nonlinear elastic body on the common boundary via the velocity and stress matching conditions. The motion of the fluid is governed by incompressible Navier-Stokes equations while the displacement of elastic structure is determined by a nonlinear elastodynamic system with boundary dissipation. The boundary dissipation is inserted in the velocity matching condition. We prove the global existence of the smooth solutions for small initial data and obtain the exponential decay of the energy of this system as well.
\\[3mm]
\\[3mm]
\end{small}
\end{quote}

\section{Introduction}
\setcounter{equation}{0}
\hskip\parindent We consider a free boundary fluid-structure interaction model which consists of the viscous incompressible fluid and the nonlinear elastic structure. Both fluid and elastic body are contained in $\Omega$ which is a smooth bounded domain in $\mathbb{R}^3$. This domain is divided into two parts by the common interface of fluid and structure where the interaction takes place. The inner part is taken over by the elastic body, denoted by $\Omega_e(t)\subset\Om$ while the fluid occupied the exterior part $\Omega_f(t)=\Omega\setminus \bar{\Omega}_e(t)$. We denote the fluid portion and the solid portion in the domain $\Omega$ at the start time by $\Omega_e$ and $\Omega_f$, respectively. Both $\Omega_e$ and $\Omega_f$ are also smooth bounded domains in $\mathbb{R}^3$. Moreover, we use the symbol $\Ga_c=\pl\Om_e\cap\pl\Om_f$ to stand for the common boundary of fluid and solid at $t=0$(For more details, see \cite{CS1,IKLT3, KT1,KT2}).
The motion of the fluid is described by the incompressible Navier-Stokes equations(see \cite{Temam 1984}):
\be u_t-\Delta u + (u\cdot\na)u+\na p =0\quad\mbox{in}\quad \Om_f(t), \label{eq1.1}\ee
\be \na\cdot u =0\quad \mbox{in}\quad \Om_f(t),\label{eq1.2} \ee
in which the vector field $u\in\mathbb{R}^3$ is the velocity of the fluid,
while the displacement of elastic structure $w\in\mathbb{R}^3$ is dominated by the following nonlinear elastodynamic system:
\begin{equation*}
w_{tt}=\div DW(Dw),
\end{equation*}
which is derived by variational methods from the action functional
\begin{equation}\label{eq1.3}
I[w]=\int^T_0\int_{\Omega_e}[\frac{1}{2}\mid w_t\mid^2-W(Dw)+|w|^2]dxd\tau,
\end{equation}
where $W(F): \mathbb{R}^{3\times 3}\rightarrow \mathbb{R}$ is the stored-energy function of the elastic material. The term $|w|^2$ in (\ref{eq1.3}) is needed for the energy estimate consideration, see \cite{IKLT2}.
The interaction occurs on the common boundary $\Gamma_c(t)$ via the natural transmission boundary conditions matching the velocity and the stress.

Fluid-structure interaction models have drawn considerable attention from both engineers and mathematical researchers. At beginning, these models were considered in a finite element framework (see(\cite{DGHL},\cite{GGCC},\cite{GGCCL}) and reference therein). Recently, the topics about the mathematical theory of existence, uniqueness and stability of solutions for such models have been becoming quite attractive. For
the linear elastic material,
local in time well-posedness of the free boundary model was first established by Coutand and Shkoller in \cite{CS1} and improved in \cite{IKLT1, KT1}  and  \cite{KT2} where there are no
dissipative  mechanisms on the interface. Then the global-in-time existence for the fluid-structure system with damping is established in \cite{IKLT2} and \cite{IKLT3} for small initial data and linear isotropic elastic material. Besides,  \cite{Qin-Yao}, the
global solutions and energy decay are also obtained for the linear wave equations with variable coefficients coupling with incompressible viscous fluid.
For other topics on fluid-linear elastic structure system, see the short review in the introduction in \cite{IKLT3}.

As for the nonlinear elastic material, in \cite{CS2}, Coutand and Shkoller developed the short-time wellposedness theory for the system in which fluid couples with some specific quasilinear elastic material.

In this paper we assume that the elastic body is a general nonlinear material to consider the global smooth solutions and energy decay of the system where the fluid interacts the nonlinear elastic body determined by the action functional $(\ref{eq1.3})$. For the fluid part of this system, we employ the method and the estimates obtained in \cite{IKLT2}. For the elastic body, we use the multiplier methods and invoke the related lemmas derived in \cite{Yao 2007}, \cite{Zhang-Yao 2008} and
\cite{Zhang-Yao 2009} to deal with our problem.

Let $\eta(x,t): \Omega\rightarrow\Omega(t)$ be the position function with $\Omega(t)=\Om$ being the different states of the system with respect to different time. With the help of position function, the incompressible Navier-Stokes equations can be reformulated in the Lagrangian framework:
\begin{eqnarray}\label{eq1.4}
&&\left\{\begin{array}{lll}
\partial_{t}v^{i}-\partial_{j}(a^{jl}a^{kl}\partial_{k}v^i)+\partial_{k}(a^{ki}q)=0 \quad &\mbox{in}\quad \Omega_f\times (0,T),\\
a^{ki}\partial_{k}v^i=0  \quad &\mbox{in}\quad \Omega_f\times (0,T),
\end{array}\right. i=1,\,2,\,3,
\end{eqnarray}
where $v(x,t)$ and $q(x,t)$ denote the Lagrangian velocity and the pressure of the fluid over the initial domain $\Omega_f$, respectively.
This means that $v(x,t)=\eta_t(x,t)=u(\eta(x,t),t)$ and $q(x,t)=p(\eta(x,t),t)$ in $\Omega_f$. The matrix $\mathbf{a}(x,t)$ is defined as the inverse of the matrix $D_x\eta(x,t)$, which means $\mathbf{a}(x,t)=(D_x\eta(x,t))^{-1}$.
Note that the Einstein summation convention is employed.
The nonlinear elastodynamic equations for the displacement function $w(x,t)=\eta(x,t)-x$ are formulated as the following:
\begin{equation}\label{eq1.5}
w_{tt}-\div DW(Dw+\mathbb{I})+w=0\quad\mbox{in}\quad \Omega_e\times (0,T),\quad i=1,2,3.
\end{equation}
We seek a solution $(v,w,q,\mathbf{a},\eta)$ to the  system \eqref{eq1.4} and \eqref{eq1.5}, where the matrix $\mathbf{a}=(a^{ij})(i,j=1,2,3)$ and $\eta\mid_{\Omega_f}$ are determined in the following way:
\begin{eqnarray}
&&\mathbf{a}_t=-\mathbf{a}:Dv:\mathbf{a},\quad\mbox{in}\quad \Omega_f\times(0,T), \label{eq1.6}\\
&&\eta_t=v \quad\mbox{in}\quad\Omega_f\times(0,T),\label{eq1.7}\\
&&\mathbf{a}(x,0)=\mathbb{I}, \quad\eta(x,0)=x,\quad\mbox{in}\quad \Omega_f,  \label{eq1.8}
\end{eqnarray}
where the symbol $``:"$ stands for the usual multiplication between matrices and $\mathbb{I}$ represents the identity matrix in
$\mathbb{R}^{3\times3}$.

Making use of the notation $\mathbf{a}$, we rewrite (\ref{eq1.4}) as
\begin{eqnarray*}
&&\left\{\begin{array}{lll}
\partial_{t}v-\div(\mathbf{a}:\mathbf{a}^{\mathbf{T}}:Dv)+\div(\mathbf{a}q)=0 \quad &\mbox{in}\quad \Omega_f\times (0,T),\\
tr(\mathbf{a}:Dv)=0  \quad &\mbox{in}\quad \Omega_f\times (0,T).\\
\end{array}\right.
\end{eqnarray*}
On the interface $\Gamma_c$ between $\Om_f$ and $\Om_e,$ we assume the transmission boundary condition
\begin{equation}\label{eq1.9}
w_t=v-\gamma w_{\nu_{\mathcal{N}}}\quad\mbox{on}\quad \Gamma_c\times(0,T),
\end{equation}
where the constant $\gamma>0$ and
$$w_{\nu_{\mathcal{N}}}=DW(Dw+\mathbb{I})\nu,$$
and the matching of stress
\begin{equation}\label{eq1.10}
w_{\nu_{\mathcal{N}}}=(\mathbf{a}:\mathbf{a}^{\mathbf{T}}:Dv)\nu-q(\mathbf{a}\nu)\quad\mbox{on}\quad\Gamma_c\times(0,T),
\end{equation} where $\nu=(\nu_1,\nu_2,\nu_3)$ is the unit outward normal with respect to $\Om_e.$
On the outside fluid boundary $\Gamma_f=\pl\Om$, we impose the non-slip condition
\begin{equation}\label{eq1.11}
v=0\quad\mbox{on}\quad\Gamma_f\times(0,T).
\end{equation}
We supplement the system (\ref{eq1.4}) and (\ref{eq1.5})  with the initial data $v(x, 0)=v_0(x)$ and $(w(x, 0),w_t(x, 0))=(w_0(x),w_1(x))$
in $\Om_f$ and $\Om_e,$ respectively.
Let
$\mathbf{H}=\{u\in L^2(\Om_f): \div u=0,\quad u\cdot\nu|_{\Gamma_f}=0\}$ and
$\mathbf{V}=\{u\in H^1(\Om_f): \div u=0,\quad u|_{\Gamma_f}=0\}.$
Based on the initial data $v_0$, the initial pressure $q_0$ is determined by solving the  problem
\begin{eqnarray}\label{eq1.12}
&&\left\{\begin{array}{lll}
\triangle q_0=-\partial_{i}v^{k}_{0}\partial_{k}v^{i}_{0}\quad &\mbox{in}\quad \Om_f,\\
Dq_{0}\cdot\nu=\triangle v_{0}\cdot\nu  \quad &\mbox{on}\quad \Gamma_f, \\
-q_0=-\partial_{j}v^{i}_{0}\nu_{j}\nu_{i}+w^{i}_{\nu_{\mathcal{N}}}\nu_i,\quad &\mbox{on}\quad \Gamma_c. \\
\end{array}\right.
\end{eqnarray}
Let the initial data $v_0\in\mathbf{V}\cap H^5(\Om_f),$ $w_0\in H^4(\Om_e),$ and $w_1\in H^3(\Om_e)$ be provided. Moreover, $v_0, w_0$
and $w_1$ are supposed to satisfy suitable compatibility conditions(for the details of the compatibility conditions, see \cite{CS2}).

We now consider the hypothesis on the stored-energy function.

{\bf (H1)}\,\,The system is in equilibrium with $w=0,$ i.e.
\begin{equation}\label{eq1.13}
DW(\mathbb{I})=0.
\end{equation}

By this assumption, we rewrite $(\ref{eq1.5})$ as following:
\be\label{eq1.14}
w_{tt}-\mathcal{N}_{w}(t)w+w=0,
\ee
where the operator $\mathcal{N}_{w}(t)\phi$ is defined for $\phi\in H^2(\Om_e)$ as
\begin{equation*}
\mathcal{N}_{w}(t)\phi=\div l_{D\phi}N_w(t),
\end{equation*}
and $N_w(t)=\int_0^1 D^2W(\mathbb{I}+sDw)ds$. And in above definition, $l_{X}K$ is defined as a $k-1$ order tensor field in $\mathbb{R}^{3\times 3}$ by
\begin{equation*}
l_{X}K(X_1,...,X_{k-1})=K(X_1,...,X_{k-1},X),\quad\mbox{for}\quad X_i,X\in\mathcal{X}(\mathbb{R}^{3\times 3}),\quad 1\leqslant i\leqslant k-1.
\end{equation*}

We need to impose the strong ellipticity condition at the zero equilibrium for the material function $W$.

{\bf (H2)}\,\,There exists a positive constant $\mu>0$ such that
\be\label{eq1.15}
D^{2}W(\mathbb{I})(b_1\otimes b_2,b_1\otimes b_2)\geqslant\mu\mid b_1\mid^2\mid b_2\mid^2,\quad\mbox{for all}\quad b_1,b_2\in\mathbb{R}^3.
\ee
The fourth order tensor $D^{2}W(F)=(\frac{\pl^2W}{\pl F_{ij}\pl F_{kl}}(F))$ is the Hessian of $W$ with respect to the variable $F=(F_{ij})\in
\mathbb{R}^{3\times 3}$.

We assume that short time solutions to problem (\ref{eq1.4}) and (\ref{eq1.5}) with the boundary damping (\ref{eq1.9})-(\ref{eq1.11}) exist if the initial data satisfy suitable compatibility conditions. This can be proved by applying the method  in \cite{CS2}.

Our main results are given as follows.

\begin{thm}\label{th1.1}
Let the assumptions ${\bf (H1)}$ and ${\bf (H2)}$ hold and $\Om_e$ be star-shaped with respect to a fixed point $x_0\in\Om_e$. Suppose that initial data $v_0\in\mathbf{V}\cap H^5(\Om_f),$ $w_0\in H^4(\Om_e),$ and $w_1\in H^3(\Om_e)$ are given small such that
the corresponding compatibility conditions hold and $\gamma\geqslant 2C>0,$ where the constant $C$ depends on the value of initial data.
Then, there exists a unique global smooth solution $(v,w,q,\mathbf{a},\eta)$ such that
\begin{eqnarray*}
&&v\in L^{\infty}([0,\infty);H^4(\Omega_f));\quad v_t\in L^{\infty}([0,\infty);H^3(\Omega_f));\\
&&v_{tt}\in L^{\infty}([0,\infty);H^2(\Omega_f)); \\
&&v_{ttt}\in L^{\infty}([0,\infty);L^2(\Om_f)),\quad Dv_{ttt}\in L^{2}([0,\infty);L^{2}(\Om_f));\\
&&\partial_{t}^{k}w\in\mathcal{C}([0,\infty);H^{4-k}(\Omega_e)),\quad k=0,1,2,3,4,
\end{eqnarray*}
with $\partial_{t}^{j}q\in L^{\infty}([0,\infty);H^{3-j}(\Omega_f))(j=0,1,2);$
$\mathbf{a},\mathbf{a}_t\in L^{\infty}([0,\infty);H^3(\Omega_f)),$ $\mathbf{a}_{tt}\in L^{\infty}([0,\infty);H^2(\Omega_f)),$ $\mathbf{a}_{ttt}\in L^{\infty}([0,\infty);L^2(\Omega_f)),$ $D\mathbf{a}_{ttt}\in L^{2}([0,\infty);L^2(\Omega_f)),$ $\pl^4_t\mathbf{a}\in L^{2}([0,\infty);L^{2}(\Omega_f))$ and $\eta|_{\Om_f}\in\mathcal{C}([0,\infty);H^{4}(\Omega_f))$.
\end{thm}

\begin{rem}
 Given initial data small the total energy of the system decays exponentially, where the total energy is defined in $(\ref{total})$ later.
\end{rem}

The proof of Theorem \ref{th1.1} will be given in Section 4.

\section{Preliminaries}
\setcounter{equation}{0}
\hskip\parindent We list some lemmas and estimates in the literature which are needed in the proof of Theorem \ref{th1.1}. Constant $C$ may be different from line to line throughout this note.

\begin{lem} \label{lem2.1}$\cite[Lemma\,\, 3.1]{IKLT1}$
Assume that $\parallel\nabla v\parallel_{L^{\infty}([0,T];H^3(\Om_f))}\leqslant M$. Let $p\in[1,\infty]$ and $1\leq i,\,j,\,k,\,l\leq 3.$
With $T\in[0,\frac{1}{CM}]$ where $C>0$ is large enough, the following statements hold:

$(i)$ $\parallel\nabla\eta\parallel_{H^3(\Omega_f)}\leqslant C\quad\mbox{for}\quad t\in[0,T];$

$(ii)$ $\parallel\mathbf{a}\parallel_{H^3(\Omega_f)}\leqslant C\quad\mbox{for}\quad t\in[0,T];$

$(iii)$ $\parallel\mathbf{a}_t\parallel_{L^p(\Omega_f)}\leqslant C\parallel\nabla v\parallel_{L^{p}(\Omega_f)}\quad\mbox{for}\quad t\in[0,T];$

$(iv)$ $\parallel\partial_{i}\mathbf{a}_t\parallel_{L^p(\Omega_f)}\leqslant C\parallel\nabla v\parallel_{L^{p_1}(\Omega_f)}\parallel\partial_{i}\mathbf{a}\parallel_{L^{p_2}(\Omega_f)}
+C\parallel\nabla\partial_{i}v\parallel_{L^p(\Omega_f)}, \quad t\in[0,T]$ where $1\leqslant p,p_1,p_2\leqslant\infty$
are given such that $\frac{1}{p}=\frac{1}{p_1}+\frac{1}{p_2}$;

$(v)$ $\parallel\partial_{ij}\mathbf{a}_t\parallel_{L^2(\Omega_f)}\leqslant C\parallel\nabla v\parallel_{H^1(\Omega_f))}^{\frac{1}{2}}\parallel\nabla v\parallel_{H^2(\Omega_f))}^{\frac{1}{2}}+C\parallel\nabla v\parallel_{H^2(\Omega_f)},\quad t\in[0,T];$

$(vi)$ $\parallel\mathbf{a}_{tt}\parallel_{L^2(\Omega_f)}\leqslant C\parallel\nabla v\parallel_{L^2(\Omega_f)}\parallel\nabla v\parallel_{L^{\infty}(\Omega_f)}+C\parallel\nabla v_t\parallel_{L^2(\Omega_f)} $ and    \\
$\parallel\mathbf{a}_{tt}\parallel_{L^3(\Omega_f)}\leqslant C\parallel v\parallel_{H^2(\Omega_f)}^2+C\parallel\nabla v_t\parallel_{L^3(\Omega_f)}$
for $t\in[0,T]$;

$(vii)$ $\parallel\mathbf{a}_{ttt}\parallel_{L^2(\Omega_f)}\leqslant C\parallel\nabla v\parallel_{H^1(\Omega_f)}^3+C\parallel\nabla v_t\parallel_{L^2(\Omega_f)}\parallel\nabla v\parallel_{L^{\infty}(\Omega_f)}+C\parallel\nabla v_{tt}\parallel_{L^2(\Omega_f),} t\in[0,T];$

$(viii)$ for every $\epsilon\in(0,\frac{1}{2}]$ and all $t\leqslant T^{*}=\min\{\frac{\epsilon}{CM^2},T\},$ we have
\begin{equation}\label{eq2.1}
\parallel\delta_{jk}-a^{jl}a^{kl}\parallel^2_{H^3(\Omega_f)}\leqslant \epsilon,\quad
\end{equation}
and
\begin{equation}\label{eq2.2}
\parallel\delta_{jk}-a^{jk}\parallel^2_{H^3(\Omega_f)}\leqslant \epsilon,\quad
\end{equation}
In particular, the form $a^{jl}a^{kl}\xi_{ij}\xi_{ik}$ satisfies the ellipticity estimates
\begin{equation}\label{eq2.3}
a^{jl}a^{kl}\xi_{ij}\xi_{ik}\geqslant \frac{1}{C}|\xi|^2,\quad \xi\in\mathbb{R}^{3\times3},
\end{equation}
for all $t\in[0,T^{*}]$ and $x\in\Omega_f$, provided $\epsilon\leqslant\frac{1}{C}$ with $C$ sufficiently large.
\end{lem}

Using a similar method carried out in Lemma \ref{lem2.1}, we have the following  estimates on $\mathbf{a}$.
As for the proof, we leave it out.

\begin{cor}\label{cor2.1}
Under the assumptions demanded in Lemma \ref{lem2.1}, it follows that for $t\in[0,T]$ with some $T>0$ and
$1\leqslant i,j,k\leqslant 3,$

$(1)$ $\|\pl_i\mathbf{a}_{tt}\|_{L^2}\leqslant C\|v\|^2_{H^2}+C\|v\|^{\frac{3}{2}}_{H^2}\|v\|^{\frac{1}{2}}_{H^3}+C\|v_t\|_{H^2};$

$(2)$ $\|\pl_t^4\mathbf{a}\|_{L^2}\leqslant C\|v\|^{\frac{7}{2}}_{H^2}\|v\|^{\frac{1}{2}}_{H^3}+C\|v\|^{\frac{1}{2}}_{H^2}\|v\|^{\frac{1}{2}}_{H^3}\|v_{tt}\|_{H^1}
+C\|v_t\|^2_{H^2}+C\|Dv_{ttt}\|_{L^2};$

$(3)$ $\|\pl_{ijk}\mathbf{a}_{t}\|_{L^2}\leqslant C\|Dv\|_{H^2}+C\|Dv\|^{\frac{1}{2}}_{H^1}\|Dv\|^{\frac{1}{2}}_{H^2}+C\|Dv\|_{H^3};$

$(4)$ $\|\pl_{ij}\mathbf{a}_{tt}\|_{L^2}\leqslant C\|v\|^2_{H^3}+C\|v\|_{H^2}+C\|v_t\|_{H^3};$

$(5)$ $\|\pl_i\mathbf{a}_{ttt}\|_{L^2}\leqslant C\|v\|^2_{H^2}+C\|v\|_{H^3}\|v_t\|_{H^1}+C\|v_t\|^2_{H^2}+C\|v\|_{H^3}\|v_{tt}\|_{H^1}+C\|v_{tt}\|_{H^1}+C\|Dv_{ttt}\|_{L^2}.$
\end{cor}

From \cite{IKLT1}, we recall  the pointwise a-priori estimates for variable coefficient Stokes system.
\begin{lem} \label{lem2.2}$\cite[Lemma\, 3.2]{IKLT1}$
Assume that $v$ and $q$ are solutions of the system
\begin{eqnarray*}
&&\left\{\begin{array}{lll}
\partial_{t}v^{i}-\partial_{j}(a^{jl}a^{kl}\partial_{k}v^i)+\partial_{k}(a^{ki}q)=0 \quad &\mbox{in}\quad \Omega_f\times (0,T),\\
a^{ki}\partial_{k}v^i=0  \quad &\mbox{in}\quad \Omega_f\times (0,T), \quad i=1,2,3\\
v^i=0,\quad\mbox{on}\quad \Gamma_f,\\
a^{jl}a^{kl}\partial_{k}v^{i}\nu_j-a^{ki}q\nu_k=(\omega^i)_{\nu_{\mathcal{N}}},\quad\mbox{on}\quad \Gamma_c,
\end{array}\right.
\end{eqnarray*}
where $a^{ji}\in L^{\infty}(\Omega_f)$ satisfies Lemma $\ref{lem2.1}$ with $\epsilon=\frac{1}{C}$ sufficiently small. Then we have
\begin{equation}\label{eq2.4}
\parallel v\parallel_{H^{s+2}(\Omega_f)}+\parallel q\parallel_{H^{s+1}(\Omega_f)}\leqslant C
\parallel v_t\parallel_{H^{s}(\Omega_f)}+C\parallel w_{\nu_{\mathcal{N}}}\parallel_{H^{s+\frac{1}{2}}(\Gamma_c)},\quad s=0,1,2
\end{equation} for $t\in(0,T).$
Moreover, we  obtain for $t\in(0,T),$ $(1)$
\begin{eqnarray}\label{eq2.5}
&&\parallel v_t\parallel_{H^3(\Omega_f)}+\parallel q_t\parallel_{H^2(\Omega_f)}\leqslant C\parallel v_{tt}\parallel_{H^1(\Omega_f)}+
C\parallel (w_t)_{\nu_{\mathcal{B}}}\parallel_{H^\frac{3}{2}(\Gamma_c)}  \nonumber \\
&&+C(\parallel v\parallel_{H^2(\Omega_f)}^{\frac{1}{2}}\parallel v\parallel_{H^3(\Omega_f)}^{\frac{1}{2}}+\|v\|_{H^3(\Omega_f)})(\parallel v\parallel_{H^3(\Omega_f)}+\parallel q\parallel_{H^2(\Omega_f)}),
\end{eqnarray}
where $T\leqslant\frac{1}{CM}$ and $C$ is sufficiently large;
and $(2)$
\begin{eqnarray*}
&&\| v_{tt}\|_{H^2(\Omega_f)}+\|q_{tt}\|_{H^1(\Omega_f)}\leqslant C\|v_{ttt}\|_{L^2(\Omega_f)}
+C\|(w_{tt})_{\nu_{\mathcal{B}}}\|_{H^{\frac{1}{2}}(\Gamma_c)}+CL(t) \nonumber \\
&&+C\|v\|_{H^2(\Om_f)}(\parallel v_t\parallel_{H^3(\Omega_f)}+\parallel q_t\parallel_{H^2(\Omega_f)}) \nonumber \\
&&+C(\parallel v\parallel_{H^3(\Omega_f)}+\parallel q\parallel_{H^2(\Omega_f)})[\|v_t\|_{H^2(\Omega_f)}
+\|v\|^2_{H^2(\Omega_f)}
+\|v\|^{\frac{3}{2}}_{H^2(\Omega_f)}\|v\|^{\frac{1}{2}}_{H^3(\Omega_f)} \nonumber \\
&&+\|v\|_{H^2(\Omega_f)}(\|v\|_{H^3(\Omega_f)}+
\|v\|^{\frac{1}{2}}_{H^3(\Omega_f)}\|v\|^{\frac{1}{2}}_{H^2(\Omega_f)})].
\end{eqnarray*}
\end{lem}

From now on, for simplicity,  we omit specifying the domains $\Omega_f$ and $\Omega_e$ in the norms involving the velocity $v$ and the displacement $w$. But we still emphasize the boundary domains $\Gamma_c$ and $\Gamma_f$.

Now we turn to the nonlinear elastodynamic system. Let $w$ be a solution of (\ref{eq1.5}) on $[0,T]$ for some $T>0$.
Set
\begin{equation*}
C(t)=DW(Dw+\mathbb{I}).
\end{equation*}
We differentiate $C(t)$ in time and have
\begin{equation}\label{eq2.6}
\dot{C}(t)=(\sum_{i,j=1}^3\frac{\pl^2W}{\pl F_{ij}\pl F_{kl}}(Dw+\mathbb{I})\frac{\pl\dot{w}^i}{\pl x^j})_{kl}=l_{Dw_t}D^2W.
\end{equation}
Then, the derivatives of order $j$ of $C(t)$ for $j=2,3,4$ with respect to $t$ are listed as follows:
\begin{eqnarray}
&&\ddot{C}(t)=l_{Dw_{tt}}D^2W+l_{Dw_t}(l_{Dw_t}D^3W),\label{eq2.7} \\
&&C^{(3)}(t)=l_{Dw_{ttt}}D^2W+3l_{Dw_t}(l_{Dw_{tt}}D^3W)+l_{Dw_t}l_{Dw_t}l_{Dw_t}D^4W, \label{eq2.8}\\
&&C^{(4)}(t)=l_{Dw_{tttt}}D^2W+4l_{Dw_t}(l_{Dw_{ttt}}D^3W)+3l_{Dw_{tt}}l_{Dw_{tt}}D^3W \nonumber  \\
&&+6l_{Dw_t}l_{Dw_t}l_{Dw_{tt}}D^4W
+l_{Dw_t}l_{Dw_t}l_{Dw_t}l_{Dw_t}D^5W. \label{eq2.9}
\end{eqnarray}

Define
\be\label{eq2.10}
\mathcal{B}_{w}(t)\phi=\div l_{D\phi}[D^2W(Dw+\mathbb{I})],\quad\mbox{for}\quad \phi\in H^2(\Om_e).
\ee
Let
\begin{equation}\label{eq2.11}
\phi_{\nu_{\mathcal{B}}}=(l_{D\phi}D^2W)\nu,\quad \mbox{for}\quad \phi\in H^2(\Om_e)
\end{equation}
be a vector field on the interface $\Ga_c$, given by the formulas
\begin{equation*}
\langle \phi_{\nu_{\mathcal{B}}},X\rangle=D^2W(D\phi,X\otimes\nu),\quad \mbox{for}\quad \phi\in H^2(\Om_e),\quad X\in\mathcal{X}(\Ga_c).
\end{equation*}

With the help of the above notations we introduce, we differentiate (\ref{eq1.5}) or (\ref{eq1.14}) in time for $j(\geqslant 1)$ times and obtain
\be\label{eq2.12}
w^{(j)}_{tt}-\mathcal{B}_{w}(t)w^{(j)}+w^{(j)}-r_{j-1}(t)=0,
\ee
where the remainder term caused by the nonlinearity of $W$ is
\be\label{eq2.13}
r_{j-1}(t)=\div(C^{(j)}(t)-\mathcal{B}_{w}(t)w^{(j)})\quad\mbox{as}\quad j\geqslant 2;\quad r_0(t)=0.
\ee
Therefore, with above preparations, we introduce the following various types of energy:

the energy for (\ref{eq1.5}) or (\ref{eq1.14})
\begin{eqnarray*}
V_0^e(t)=\frac{1}{2}(\|w_t\|^2_{L^2}+\|w\|^2_{L^2}+\int_{\Om_e}\langle l_{Dw}N_w(t),Dw\rangle_{\mathbb{R}^{3\times3}}dx),
\end{eqnarray*}
the first level energy of the whole system
\begin{eqnarray*}
V_0(t)=V_0^e(t)+\frac{1}{2}\|v\|^2_{L^2},
\end{eqnarray*}
for $j\geqslant 1$, the energy for (\ref{eq2.12})
\begin{eqnarray*}
V_j^e(t)=\frac{1}{2}(\|w_t^{(j)}\|^2_{L^2}+\|w^{(j)}\|^2_{L^2}+\int_{\Om_e}\langle l_{Dw^{(j)}}D^2W,Dw^{(j)}\rangle_{\mathbb{R}^{3\times3}}dx),
\end{eqnarray*}
the j-th level energy of the fluid-structure interaction model
\begin{eqnarray*}
V_j(t)=V_j^e(t)+\frac{1}{2}\|v^{(j)}\|^2_{L^2},
\end{eqnarray*}
and the total energy of the whole system and the elastic body
\begin{eqnarray*}
\mathcal{Q}(t)=\sum_{j=0}^3V_j(t),\quad\mbox{and}\quad \mathcal{E}^e(t)=\sum_{k=0}^4\|w^{(k)}(t)\|^2_{H^{4-k}(\Om_e)}.
\end{eqnarray*}
Moreover, we denote a remainder term which will appear by $L(t)=\sum_{k=3}^{8}(\mathcal{E}^e(t))^{\frac{k}{2}}$.

Now, with necessary notation introduced above, we give some lemmas which will be frequently utilized.

Suppose that the assumption {\bf (H2)} holds. Hence, there exists a constant $\gamma_0>0$ such that if $w\in H^1(\Om_e)$ with
\be\label{eq2.14}
\sup_{x\in\overline{\Om}_e}\mid Dw\mid\leqslant \gamma_0,
\ee
then
\be\label{eq2.15}
D^2W(Dw+\mathbb{I})(b_1\otimes b_2,b_1\otimes b_2)\geqslant\mu\mid b_1\mid^2\mid b_2\mid^2,
\ee
for all $b_1,b_2\in\mathbb{R}^3$. Thus, it yields the following lemma, by G\"{a}rding's inequality.
\begin{lem}\label{lem2.3}$\cite[Lemma\, 2.1]{Zhang-Yao 2009}$
Let {\bf (H2)} hold. Then there is $\gamma_0>0$ such that, if $w\in H^1(\Om_e)$ is such that the condition (\ref{eq2.14}) holds, then
\be\label{eq2.16}
\int_{\Om_e}D^2W(Dw+\mathbb{I})(D\phi,D\phi)dx\geqslant c_{\gamma_0}\|D\phi\|^2_{L^2},\quad \phi\in H^1(\Om_e),
\ee
for some $c_{\gamma_0}>0$.
\end{lem}

By slightly modifying the proofs in \cite{Yao 2007, Zhang-Yao 2008}, we obtain the following lemmas.
Before the statement of these lemmas, we collect a few basic properties of Sobolev space which we'll use often first.

$(i)$ Let $s_1>s_2\geqslant 0$ and $\mathcal{D}$ be a bounded, open set with smooth boundary in $\mathbb{R}^n$($n\geqslant 2$). For any $\epsilon>0$ there is $c_{\epsilon}>0$ such that
\be\label{eq2.17}
\|\phi\|^2_{H^{s_2}(\mathcal{D})}\leqslant \epsilon \|\phi\|^2_{H^{s_1}(\mathcal{D})}+c_{\epsilon}\|\phi\|^2_{L^{2}(\mathcal{D})},\quad\forall \phi\in H^{s_1}(\mathcal{D}).
\ee

$(ii)$ If $s>\frac{n}{2}$, then for each $k=0,...,$ we have $H^{s+k}(\mathcal{D})\subset\mathcal{C}^k(\overline{\mathcal{D}})$ with continuous inclusion.

$(iii)$ If $r\triangleq \min\{s_1,s_2,s_1+s_2-[\frac{n}{2}]-1\}\geqslant 0$, then there is a constant $c>0$ such that
\be\label{eq2.18}
\|f_1f_2\|_{H^r(\mathcal{D})}\leqslant c\|f_1\|_{H^{s_1}(\mathcal{D})}\|f_2\|_{H^{s_2}(\mathcal{D})}\quad \forall f_1\in H^{s_1}(\mathcal{D}),f_2\in H^{s_2}(\mathcal{D}).
\ee

$(iv)$ Let $s_j\geqslant 0, j=1,...,k$ and $r\triangleq \min_{1\leqslant i\leqslant k}\min_{j_1\leqslant...\leqslant j_i}\{s_{j_1}+\cdot\cdot\cdot+s_{j_i}-(i-1)\times([\frac{n}{2}]+1)\}\geqslant 0.$ Then there is a constant $c>0$ such that
\be\label{eq2.19}
\|f_1\cdot\cdot\cdot f_k\|_{H^r(\mathcal{D})}\leqslant c\|f_1\|_{H^{s_1}(\mathcal{D})}\cdot\cdot\cdot\|f_k\|_{H^{s_k}(\mathcal{D})}\quad \forall f_j\in H^{s_j}(\mathcal{D}),1\leqslant j\leqslant k.
\ee

\begin{lem}\label{lem2.4}$\cite[Lemma\, 2.1]{Yao 2007}$
Let $\gamma_0>0$ be given and $\phi\in H^1(\mathcal{D})$ be such that the condition (\ref{eq2.14}) hold. Let $f(\cdot,\cdot)$ be a smooth function on $\overline{\mathcal{D}}\times\mathbb{R}^n$ and $F(x)=f(x,D\phi).$ Then there is $c_{\gamma_0}>0,$ depending on
$\gamma_0,$ such that
\be\label{eq2.20}
\|F\|_{H^k(\mathcal{D})}\leqslant c_{\gamma_0}\sum_{j=1}^k(1+\|\phi\|_{H^m(\mathcal{D})})^j,
\ee
for $0\leqslant k\leqslant m-1$ and $m\geqslant[\frac{n}{2}]+3.$ Moreover, if $f(x,0)=0$ for $x\in\mathcal{D}$ and $\|\phi\|_{H^m(\mathcal{D})}\leqslant \gamma_0,$ then there is $c_{\gamma_0}>0$ such that
\be\label{eq2.21}
\|F\|_{H^k(\mathcal{D})}\leqslant c_{\gamma_0}\|\phi\|_{H^{k+1}(\mathcal{D})},\quad 0\leqslant k\leqslant m-1.
\ee
\end{lem}

\begin{lem}\label{lem2.5}$\cite[Lemma\, 2.2]{Yao 2007}$
Let $\gamma_0>0$ be given and $w$ satisfy the problem (\ref{eq1.5}) on the interval $[0,T]$ for some $T>0$ such that
\be\label{eq2.22}
\sup_{0\leqslant t\leqslant T}\|w(t)\|_{H^m(\Om_e)}\leqslant \gamma_0.
\ee
For $1\leqslant j\leqslant 2,$ let
\begin{equation*}
r_{j}(t)=\div(C^{(j+1)}(t)-\mathcal{B}_{w}(t)w^{(j+1)}),\quad r_{j,\Ga_c}=[C^{(j+1)}(t)-\mathcal{B}_{w}(t)w^{(j+1)}]\nu,
\end{equation*}
where the operator $\mathcal{B}_{w}(t)$ is given by (\ref{eq2.10}).
Then,
\be\label{eq2.23}
\|r_{j}(t)\|^2_{H^{2-j}(\Om_e)},\|r_{j,\Ga_c}(t)\|^2_{H^{\frac{5}{2}-j}(\Ga_c)}\leqslant c_{\gamma_0}\sum_{k=2}^3\mathcal{E}^k(t),
\quad j=1,2.
\ee
\end{lem}

\begin{lem}\label{lem2.6}$\cite[Lemma\, 2.3]{Yao 2007}$
Let $\gamma_0>0$ be given and $w$ satisfy the problem (\ref{eq1.5}) on the interval $[0,T]$ for some $T>0$ such that
(\ref{eq2.22}) is true. Then there is $c_{\gamma_0}>0$, which depends on $\gamma_0$, such that
\be\label{eq2.24}
\|\phi\|^2_{H^{k+1}(\Om_e)}\leqslant c_{\gamma_0}(\|\mathcal{B}_{w}(t)\phi\|^2_{H^{k-1}(\Om_e)}+\|\phi_{\nu_{\mathcal{B}}}\|^2_{H^{k-\frac{1}{2}}(\Ga_c)}+\|\phi\|^2_{H^{k}(\Om_e)}),
\quad 1\leqslant k\leqslant 3,
\ee
for $\phi\in H^{k+1}(\Om_e)$ and $t\in[0,T]$.
\end{lem}

Combining all the lemmas from (\ref{lem2.4}) to (\ref{lem2.6}), we arrive at the following theorem for $\mathcal{Q}^e(t)=\sum_{j=0}^3V_j^e(t)$ and $\mathcal{E}^e(t)$.

\begin{thm}\label{th2.1}$\cite[Theorem\, 2.1]{Yao 2007}$
Let $\gamma_0>0$ be given and $w$ satisfy the problem (\ref{eq1.5}) on the interval $[0,T]$ for some $T>0$ such that
(\ref{eq2.22}) holds. Then there are constants $c_{0,\gamma_0}>0$ and $c_{\gamma_0}>0$, which depend on $\gamma_0$, such that
\be\label{eq2.25}
c_{0,\gamma_0}\mathcal{Q}^e(t)\leqslant\mathcal{E}^e(t)\leqslant c_{\gamma_0}\mathcal{Q}^e(t)+c_{\gamma_0}L(t)
+c_{\gamma_0}\sum_{j=0}^2\|w^{(j)}_{\nu_{\mathcal{B}}}\|^2_{H^{\frac{5}{2}-j}(\Ga_c)},
\ee
\be\label{eq2.26}
-c_{\gamma_0}L(t)-c_{\gamma_0}L^{\frac{2}{3}}(t)+\sum_{j=0}^3\int_{\Ga_c}\langle (l_{Dw^{(j)}}D^2W)\cdot\nu,w^{(j+1)}\rangle d\sigma\leqslant\dot{\mathcal{Q}}^e(t),
\ee
\be\label{eq2.27}
\dot{\mathcal{Q}}^e(t)
\leqslant c_{\gamma_0}L(t)+c_{\gamma_0}L^{\frac{2}{3}}(t)+\sum_{j=0}^3\int_{\Ga_c}\langle (l_{Dw^{(j)}}D^2W)\cdot\nu,w^{(j+1)}\rangle d\sigma,
\ee
and
\be\label{eq2.28}
V_0^e(t)\leqslant c_{\gamma_0}[\int^t_0L(\tau)+\|w_{\nu_{\mathcal{N}}}\|^{\frac{3}{2}}_{L^2(\Ga_c)}d\tau
+V_0^e(0)]e^t.
\ee
\end{thm}

Thus, all the preparations for the proof of Theorem \ref{th1.1} have been made.

\section{A-priori estimates}
\setcounter{equation}{0}
\hskip\parindent
We derive a priori estimates for the global existence of  solutions to the dissipative fluid-nonlinear structure system when the initial data are sufficiently small.

Assume that
\begin{equation*}
\parallel v_0\parallel^2_{H^1}, \parallel v_t(0)\parallel^2_{H^1}, \parallel v_{tt}(0)\parallel^2_{H^1},\|v_{ttt}(0)\|^2_{L^2},
\parallel w_0\parallel^2_{H^4},\parallel w_1\parallel^2_{H^3}\leqslant \epl,
\end{equation*}
where $\epsilon>0$ is given small enough,i.e. $\epl<\gamma_0.$

We demand several auxiliary estimates involving different levels of energy.

\subsection{First level estimates}
\hskip\parindent
As defined above, $V_0(t)$ is the first level energy. We deal with this energy in this subsection.

\begin{lem}\label{lem3.1}
Suppose that all the assumptions of Theorem \ref{th1.1} are true,
the following inequality holds for $t\in[0,T]$
\begin{equation}\label{eq3.1}
V_0(t)+\int_0^{t}D_0(\tau)d\tau\leqslant V_0(0)+C_{\gamma_0}\int_0^tL(\tau)d\tau,
\end{equation}
where
\begin{equation}\label{eq3.2}
D_0(t)=\frac{1}{C}\parallel Dv(t)\parallel^2_{L^2}
+\gamma\parallel w_{\nu_{\mathcal{N}}}(t)\parallel^2_{L^2(\Gamma_c)}
\end{equation}
is a dissipative term.
\end{lem}

{\bf Proof} \,\,\, Take the $L^2$-inner product of (\ref{eq1.4}) with $v$ and (\ref{eq1.5}) with $w_t$, respectively. Hence, by (\ref{eq2.3}), we have
\be\label{eq3.3}
\frac{1}{2}\frac{d}{dt}\|v\|^2_{L^2}+\frac{1}{C}\|Dv\|^2_{L^2}\leqslant-\int_{\Ga_c}\langle(\mathbf{a}:\mathbf{a}^{\mathbf{T}}:Dv)\nu
-q(\mathbf{a}\nu),v\rangle d\sigma
\ee
and
\beq\label{eq3.4}
&&\frac{1}{2}\frac{d}{dt}(\|w_t\|^2_{L^2}+\|w\|^2_{L^2}+\int_{\Om_e}\langle l_{Dw}N_w(t),Dw\rangle dx)
-\frac{1}{2}\int_{\Om_e}N'_w(t)(Dw_t,Dw,Dw)dx   \nonumber \\
&&=\int_{\Ga_c}\langle w_{\nu_{\mathcal{N}}},w_t\rangle d\sigma,
\eeq
where $N'_w(t)(Dw_t,Dw,Dw)=\int_0^1s\frac{\pl^3W}{\pl F_{ij}\pl F_{kl}\pl F_{\alpha\beta}}(Dw+\mathbb{I})ds\frac{\pl w_t^{\alpha}}{\pl x^{\beta}}\frac{\pl w^k}{\pl x^l}\frac{\pl w^i}{\pl x^j}$ as $i,j,k,l,\alpha,\beta=1,2,3.$

Note that
\begin{equation*}
\int_0^t\int_{\Om_e}N'_w(t)(Dw_t,Dw,Dw)dxd\tau\leqslant C_{\gamma_0}\int_0^tL(\tau)d\tau.
\end{equation*}

Add (\ref{eq3.3}) and (\ref{eq3.4}) together and integrate in time from $0$ to $t$. With the help of boundary condition
(\ref{eq1.9}) and (\ref{eq1.10}), we arrive at the resulted inequality.   \hfill$\Box$

In order to proceed further, we have to derive the multiplier identities which will play a central role in our calculations.

Denote $\hat{w}=w^{(j)}$ for $0\leqslant j\leqslant 3$. From (\ref{eq1.14}) and (\ref{eq2.12}), $\hat{w}$ satisfies the following
equations:
\be\label{eq3.5}
\hat{w}_{tt}-\mathcal{A}_{w}(t)\hat{w}+\hat{w}-r(t)=0,
\ee
where $r(t)=r_{j-1}(t)$ as $j=2,3;$ $r(t)=0$ when $j=0,1$ and the operator $\mathcal{A}_{w}(t)=\mathcal{B}_{w}(t)$ if $j\geqslant 1,$
or else $\mathcal{A}_{w}(t)=\mathcal{N}_{w}(t)$ if $j=0.$

\begin{lem}\label{lem3.2}
Let $\hat{w}$ be a solution of (\ref{eq3.5}) and $H$ be a vector field on $\overline{\Om}_e$ with $H(\hat{w})=D_H\hat{w}.$ And let the scalar function
$\xi(x)\in \mathcal{C}^1(\overline{\Om}_e).$ Then we have for $\varrho>0,$
\beq\label{eq3.6}
&&\int_{s}^t\int_{\Ga_c}(\mid\hat{w}_t\mid^2-A_{w}(t)(D\hat{w},D\hat{w})-\mid\hat{w}\mid^2)\langle H,\nu\rangle d\sigma d\tau
+\int_{s}^t\int_{\Ga_c}\langle \hat{w}_{\nu_{\mathcal{A}}},2H(\hat{w})+\varrho\hat{w}\rangle d\sigma d\tau \nonumber \\
&&=(\hat{w}_t,2H(\hat{w})
+\varrho\hat{w})_{L^2}\mid^t_s+\int_{s}^t\int_{\Om_e}(\mid\hat{w}_t\mid^2-A_{w}(t)(D\hat{w},D\hat{w})-\mid\hat{w}\mid^2)(\div H-\varrho)dxd\tau \nonumber \\
&&+2\int_{s}^t\int_{\Om_e}A_{w}(t)(D\hat{w},D\hat{w}:DH)dxd\tau
-\int_{s}^t\int_{\Om_e}DA_{w}(t)(D\hat{w},D\hat{w},D_HDw)dxd\tau  \nonumber \\
&&-\int_{s}^t\int_{\Om_e}\langle r(t),2H(\hat{w})+\varrho\hat{w}\rangle dxd\tau
\eeq
and
\beq\label{eq3.7}
&&(\hat{w}_t,\xi\hat{w})_{L^2}\mid^t_s-\int_{s}^t\int_{\Om_e}\langle r(t),\xi\hat{w}\rangle dxd\tau
+\int_{s}^t\int_{\Om_e}\xi(\mid\hat{w}\mid^2-\mid\hat{w}_t\mid^2)dxd\tau \nonumber \\
&&+\int_{s}^t\int_{\Om_e}\xi A_{w}(t)(D\hat{w},D\hat{w})dxd\tau+\int_{s}^t\int_{\Om_e}A_{w}(t)(D\hat{w},D\xi\otimes\hat{w})dxd\tau \nonumber \\
&&=\int_{s}^t\int_{\Ga_c}\xi\langle\hat{w}_{\nu_{\mathcal{A}}},\hat{w}\rangle d\sigma d\tau,
\eeq
where $A_{w}(t)=D^2W(Dw+\mathbb{I}), DA_{w}(t)=D^3W$ if $j\geqslant 1$ or else $A_{w}(t)=N_w(t), DA_{w}(t)=N'_w(t)$ if $j=0.$
\end{lem}
{\bf Proof} \,\,\, For (\ref{eq3.6}), we take $2H(\hat{w})+\varrho\hat{w}$ as the multiplier.
Multiply (\ref{eq3.5}) by $2H(\hat{w})+\varrho\hat{w}$ and integrate by parts over $(s,t)\times \Om_e.$
Thus, we have
\beq\label{eq3.8}
&&\int_{s}^t\int_{\Om_e}\langle\hat{w}_{tt},2H(\hat{w})+\varrho\hat{w}\rangle dxd\tau
=(\hat{w}_t,2H(\hat{w})+\varrho\hat{w})_{L^2}\mid^t_s+\int_{s}^t\int_{\Om_e}\mid\hat{w}_t\mid^2(\div H-\varrho)dxd\tau \nonumber \\
&&-\int_{s}^t\int_{\Ga_c}\mid\hat{w}_t\mid^2\langle H,\nu\rangle d\sigma d\tau;
\eeq
\beq\label{eq3.9}
&&-\int_{s}^t\int_{\Om_e}\langle\mathcal{A}_{w}(t)\hat{w},2H(\hat{w})+\varrho\hat{w}\rangle dxd\tau
=\int_{s}^t\int_{\Om_e}\langle l_{D\hat{w}}A_{w}(t),2DH(\hat{w})+\varrho D\hat{w}\rangle dxd\tau \nonumber \\
&&-\int_{s}^t\int_{\Ga_c}\langle\hat{w}_{\nu_{\mathcal{A}}},2H(\hat{w})+\varrho\hat{w}\rangle d\sigma d\tau
=2\int_{s}^t\int_{\Om_e}A_{w}(t)(D\hat{w},D\hat{w}:DH)dxd\tau \nonumber \\
&&+\int_{s}^t\int_{\Om_e}(\varrho-\div H)A_{w}(t)(D\hat{w},D\hat{w})dxd\tau
-\int_{s}^t\int_{\Om_e}DA_{w}(t)(D\hat{w},D\hat{w},D_HDw)dxd\tau  \nonumber \\
&&\int_{s}^t\int_{\Ga_c}A_{w}(t)(D\hat{w},D\hat{w})\langle H,\nu\rangle d\sigma d\tau
-\int_{s}^t\int_{\Ga_c}\langle \hat{w}_{\nu_{\mathcal{A}}},2H(\hat{w})+\varrho\hat{w}\rangle d\sigma d\tau,
\eeq
with the help of
\begin{equation*}
D[H(\hat{w})]=D\hat{w}:DH+D_HD\hat{w}
\end{equation*}
and
\begin{eqnarray*}
&&\div[A_{w}(t)(D\hat{w},D\hat{w})H]=A_{w}(t)(D\hat{w},D\hat{w})\div H
+2A_{w}(t)(D\hat{w},D_HD\hat{w})  \\
&&+DA_{w}(t)(D\hat{w},D\hat{w},D_HDw);
\end{eqnarray*}
and
\begin{eqnarray} \label{eq3.10}
\int_{s}^t\int_{\Om_e}\langle\hat{w},2H(\hat{w})+\varrho\hat{w}\rangle dxd\tau
&=&\int_{s}^t\int_{\Om_e}\varrho\mid\hat{w}\mid^2dxd\tau
-\int_{s}^t\int_{\Om_e}\div H\mid\hat{w}\mid^2dxd\tau \nonumber \\
&+&\int_{s}^t\int_{\Ga_c}\mid\hat{w}\mid^2\langle H,\nu\rangle d\sigma d\tau
\end{eqnarray}
Combining (\ref{eq3.8})-(\ref{eq3.10}), we acquire (\ref{eq3.6}).

As for (\ref{eq3.7}), we regard $\xi(x)\hat{w}$ as the multiplier. Carry out similar calculations and we'll arrive at the multiplier identity (\ref{eq3.7}). \hfill$\Box$

Taking advantage of the above multiplier identities, we establish the following Lemma. Moreover, we intend to deal with the similar
estimates in subsection 3.1-3.4 in a unified way.

The related boundary conditions of (\ref{eq3.5}) are the following:
\be\label{eq3.11}
\hat{w}_t=\hat{v}-\gamma(\hat{w}_{\nu_{\mathcal{A}}}+\hat{r}(t)),
\ee
where $\hat{v}=v^{(j)}$ for the related $0\leqslant j\leqslant 3,$
and
\be\label{eq3.12}
\hat{w}_{\nu_{\mathcal{A}}}=\pl^{(j)}_t(\mathbf{a}:\mathbf{a}^{\mathbf{T}}:Dv)\nu-\pl^{(j)}_t(q\mathbf{a})\nu-\hat{r}(t),
\ee
where $0\leqslant j\leqslant 3$ and $\hat{r}(t)=0$ if $j=0,1$ or else $\hat{r}(t)=r_{j-1,\Ga_c}=[C^{(j)}(t)-\mathcal{B}_{w}(t)w^{(j)}]\nu$
if $j=2,3.$ Besides, in order to simplify the notations, we denote $V(t)=V_j(t)$ as $0\leqslant j\leqslant 3$ and
$V^e(t)=V_j^e(t).$

\begin{lem}\label{lem3.3}
Under the assumptions of Theorem \ref{th1.1},
it holds for $0\leqslant s<t\leqslant T$ with some $T>0$ that there exists a positive constant $C$ such that
\beq\label{eq3.13}
&&\int_{s}^tV(\tau)d\tau
\leqslant CV(s)+CV(t)+C\int_{s}^tL(\tau)d\tau
+C\int_{s}^t\int_{\Ga_c}\mid\hat{w}_{\nu_{\mathcal{A}}}\mid^2d\sigma d\tau \nonumber \\
&&+C\int_{s}^t\|D\hat{v}\|^2_{L^2}d\tau.
\eeq
\end{lem}

{\bf Proof} \,\,\, Submit $H=x-x_0$ and $\varrho=3-\varepsilon$, where $0<\varepsilon<1$ is to be determined, into (\ref{eq3.6}).
Note that $\div H=3$ and $DH=\mathbb{I}$. Hence, we have
\beq\label{eq3.14}
&&\int_{s}^t\int_{\Ga_c}(\mid\hat{w}_t\mid^2-A_{w}(t)(D\hat{w},D\hat{w})-\mid\hat{w}\mid^2)\langle H,\nu\rangle d\sigma d\tau
+\int_{s}^t\int_{\Ga_c}\langle \hat{w}_{\nu_{\mathcal{A}}},2H(\hat{w})+\varrho\hat{w}\rangle d\sigma d\tau \nonumber \\
&&=(\hat{w}_t,2H(\hat{w})
+\varrho\hat{w})_{L^2}\mid^t_s+\varepsilon\int_{s}^t\int_{\Om_e}(\mid\hat{w}_t\mid^2-A_{w}(t)(D\hat{w},D\hat{w})-\mid\hat{w}\mid^2)dxd\tau \nonumber \\
&&+2\int_{s}^t\int_{\Om_e}A_{w}(t)(D\hat{w},D\hat{w})dxd\tau
-\int_{s}^t\int_{\Om_e}DA_{w}(t)(D\hat{w},D\hat{w},D_HDw)dxd\tau  \nonumber \\
&&-\int_{s}^t\int_{\Om_e}\langle r(t),2H(\hat{w})+\varrho\hat{w}\rangle dxd\tau.
\eeq

Now we estimate the terms in (\ref{eq3.14}).

By Cauchy-Schwartz inequality and (\ref{eq2.15}) or Lemma \ref{lem2.3}, we have
\be\label{eq3.15}
\mid(\hat{w}_t,2H(\hat{w})+\varrho\hat{w})_{L^2}\mid^t_s\mid\leqslant CV(s)+CV(t).
\ee

With the help of Lemma \ref{lem2.5} and Young inequality, we derive that
\be\label{eq3.16}
\mid\int_{s}^t\int_{\Om_e}\langle r(t),2H(\hat{w})+\varrho\hat{w}\rangle dxd\tau\mid\leqslant \epsilon_1\int_{s}^tV^e(\tau)d\tau
+C_{\epsilon_1}\int_{s}^tL(\tau)d\tau,
\ee
where $\epsilon_1>0$ is sufficiently small and to be determined later on.

Applying the property of Stored-function $W$ and the Sobolev inequality $\|\phi\|_{L^{\infty}}\leqslant C\|\phi\|_{H^2}$, it follows that
\be\label{eq3.17}
\mid\int_{s}^t\int_{\Om_e}DA_{w}(t)(D\hat{w},D\hat{w},D_HDw)dxd\tau\mid\leqslant C\int_{s}^tL(\tau)d\tau.
\ee

Using the Cauchy-Schwartz inequality, Young inequality and Lemma \ref{lem2.3} again, we have
\beq\label{eq3.18}
\mid\int_{s}^t\int_{\Ga_c}\langle \hat{w}_{\nu_{\mathcal{A}}},2H(\hat{w})+\varrho\hat{w}\rangle d\sigma d\tau\mid
&\leqslant& \epsilon_2\int_{s}^t\int_{\Ga_c}[A_{w}(t)(D\hat{w},D\hat{w})+\mid\hat{w}\mid^2]d\sigma d\tau \nonumber \\
&+&C_{\epsilon_2}\int_{s}^t\int_{\Ga_c}\mid\hat{w}_{\nu_{\mathcal{A}}}\mid^2d\sigma d\tau,
\eeq
in which $\epsilon_2>0$ is also to be determined.

By the star-shaped condition for the domain $\Om_e$, it implies that there exists a constant $\rho_0>0$ such that
\be\label{eq3.19}
\langle x-x_0,\nu\rangle\geqslant \rho_0.
\ee
Therefore, we take $\epsilon_2\leqslant\frac{\rho_0}{2}.$

Insert (\ref{eq3.15})-(\ref{eq3.19}) into (\ref{eq3.14}) and thus we deduce that
\beq\label{eq3.20}
&&\varepsilon\int_{s}^t\int_{\Om_e}\mid\hat{w}_t\mid^2dxd\tau
+(2-\varepsilon)\int_{s}^t\int_{\Om_e}A_{w}(t)(D\hat{w},D\hat{w})dxd\tau \nonumber \\
&&+\frac{\rho_0}{2}\int_{s}^t\int_{\Ga_c}[A_{w}(t)(D\hat{w},D\hat{w})+\mid\hat{w}\mid^2]d\sigma d\tau
\leqslant\varepsilon\int_{s}^t\int_{\Om_e}\mid\hat{w}\mid^2dxd\tau \nonumber \\
&&+\int_{s}^t\int_{\Ga_c}\mid\hat{w}_t\mid^2\langle H,\nu\rangle d\sigma d\tau
+C_{\epsilon_2}\int_{s}^t\int_{\Ga_c}\mid\hat{w}_{\nu_{\mathcal{A}}}\mid^2d\sigma d\tau \nonumber \\
&&+CV(s)+CV(t)+\epsilon_1\int_{s}^tV^e(\tau)d\tau+C_{\epsilon_1}\int_{s}^tL(\tau)d\tau.
\eeq

Note that the Poincar\'{e} inequality
\begin{equation*}
\|\phi\|^2_{L^2}\leqslant C\|D\phi\|^2_{L^2}+C\int_{\Ga_c}\mid\phi\mid^2d\sigma
\end{equation*}
and set $\varepsilon$ small enough so that $2-\varepsilon(1+2CC_{\gamma_0})\geqslant1$ and $\frac{\rho_0}{2}-2\varepsilon C\geqslant \frac{\rho_0}{4}$ hold. Then it leads to
\beq\label{eq3.21}
&&\varepsilon\int_{s}^t\int_{\Om_e}\mid\hat{w}_t\mid^2dxd\tau
+\int_{s}^t\int_{\Om_e}A_{w}(t)(D\hat{w},D\hat{w})dxd\tau+\varepsilon\int_{s}^t\int_{\Om_e}\mid\hat{w}\mid^2dxd\tau \nonumber \\
&&+\frac{\rho_0}{4}\int_{s}^t\int_{\Ga_c}[A_{w}(t)(D\hat{w},D\hat{w})+\mid\hat{w}\mid^2]d\sigma d\tau
\leqslant
\int_{s}^t\int_{\Ga_c}\mid\hat{w}_t\mid^2\langle H,\nu\rangle d\sigma d\tau \nonumber \\
&&+C_{\epsilon_2}\int_{s}^t\int_{\Ga_c}\mid\hat{w}_{\nu_{\mathcal{A}}}\mid^2d\sigma d\tau
+CV(s)+CV(t)+\epsilon_1\int_{s}^tV^e(\tau)d\tau+C_{\epsilon_1}\int_{s}^tL(\tau)d\tau.
\eeq

Let the constant $\epsilon_1<\frac{\varepsilon}{2}$. Thus we have
\beq\label{eq3.22}
&&\int_{s}^t\int_{\Om_e}\mid\hat{w}_t\mid^2dxd\tau
+\int_{s}^t\int_{\Om_e}A_{w}(t)(D\hat{w},D\hat{w})dxd\tau+\int_{s}^t\int_{\Om_e}\mid\hat{w}\mid^2dxd\tau \nonumber \\
&&+\frac{\rho_0}{4\varepsilon}\int_{s}^t\int_{\Ga_c}[A_{w}(t)(D\hat{w},D\hat{w})+\mid\hat{w}\mid^2]d\sigma d\tau
\leqslant
\frac{C}{\varepsilon}\int_{s}^t\int_{\Ga_c}\mid\hat{w}_t\mid^2d\sigma d\tau \nonumber \\
&&+\frac{C_{\epsilon_2}}{\varepsilon}\int_{s}^t\int_{\Ga_c}\mid\hat{w}_{\nu_{\mathcal{A}}}\mid^2d\sigma d\tau
+\frac{C}{\varepsilon}V(s)+\frac{C}{\varepsilon}V(t)+\frac{C_{\epsilon_1}}{\varepsilon}\int_{s}^tL(\tau)d\tau.
\eeq

Add $\frac{1}{2}\int_s^t\int_{\Om_e}\mid\hat{v}\mid^2dxd\tau$ to both sides of (\ref{eq3.22}) and submit the boundary condition (\ref{eq3.11}) into (\ref{eq3.22}).
Therefore, thanks to the Poincar\'{e} inequality $\|\phi\|_{L^2}\leqslant C\|D\phi\|_{L^2}$ for $\phi\in H^1_{\Ga_f}(\Om_f)$ and
(\ref{eq2.23}) in Lemma \ref{lem2.5},
 (\ref{eq3.13}) can be obtained directly from (\ref{eq3.22}).  \hfill$\Box$

According to Lemma \ref{lem3.3}, for the first order energy,
it implies that
\beq\label{eq3.23}
&&\int_{s}^tV_0(\tau)d\tau
\leqslant CV_0(s)+CV_0(t)+C\int_{s}^tL(\tau)d\tau
+C\int_{s}^t\int_{\Ga_c}\mid w_{\nu_{\mathcal{N}}}\mid^2d\sigma d\tau \nonumber \\
&&+C\int_{s}^t\|Dv\|^2_{L^2}d\tau.
\eeq
Next, we multiply (\ref{eq3.23}) by sufficiently small $\epsilon_3>0$ and add the resulted inequality to (\ref{eq3.1}).
Thus, we get the following lemma.
\begin{lem}\label{lem3.4}
Assume that all of the hypotheses in Theorem \ref{th1.1} hold. Then it's true for $t\in [0,T]$ with some $T>0$ that
\be\label{eq3.24}
V_0(t)+\int^t_0V_0(\tau)d\tau\leqslant CV_0(0)+C\int_{0}^tL(\tau)d\tau.
\ee
\end{lem}

\begin{rem}
If the solution of this fluid-nonlinear structure interaction system can exist for all time and all the assumptions in Theorem \ref{th1.1} hold, then we may infer from the above lemma that the first level energy $V_0(t)$ decays exponentially.
\end{rem}

\subsection{Second level estimates}
\hskip\parindent
As we have defined in Section 2, the second order energy
\begin{equation*}
V_1(t)=V_1^e(t)+\frac{1}{2}\|v_t\|^2_{L^2},\quad\mbox{with}\quad
V_1^e(t)=\frac{1}{2}(\|w_{tt}\|^2_{L^2}+\|w_t\|^2_{L^2}+\int_{\Om_e}\langle l_{Dw_t}D^2W,Dw_t\rangle dx),
\end{equation*}
and the corresponding dissipation term
\begin{equation*}
D_1(t)=\frac{1}{C}\parallel Dv_t(t)\parallel^2_{L^2}
+\gamma\parallel (w_t)_{\nu_{\mathcal{B}}}(t)\parallel^2_{L^2(\Gamma_c)}.
\end{equation*}

To start with deriving the similar estimates with (\ref{eq3.1}) and (\ref{eq3.24}) in above subsection, we differentiate the whole system in time. It follows that
\begin{eqnarray}
&&\partial_{t}v_t-\partial_{t}\div(\mathbf{a}:\mathbf{a}^{\mathbf{T}}:Dv)+\partial_{t}\div(\mathbf{a}q)=0 \quad \mbox{in}\quad \Omega_f\times (0,T),\label{eq3.25}\\
&&tr(\mathbf{a}_t:Dv)+tr(\mathbf{a}:Dv_t)=0  \quad \mbox{in}\quad \Omega_f\times (0,T),\label{eq3.26}\\
&&w_{ttt}-\mathcal{B}_{w}(t)w_t+w_t=0 \quad \mbox{in}\quad \Omega_e\times (0,T).\label{eq3.27}
\end{eqnarray}

Besides, the boundary conditions with respect to (\ref{eq3.25})-(\ref{eq3.27}) are
\beq
&&w_{tt}=v_t-\gamma(w_t)_{\nu_{\mathcal{B}}} \quad\mbox{on}\quad \Ga_c\times (0,T),\label{eq3.28}\\
&&(w_t)_{\nu_{\mathcal{B}}}=\pl_t(\mathbf{a}:\mathbf{a}^{\mathbf{T}}:Dv)\nu-\pl_t(q\mathbf{a})\nu \quad\mbox{on}\quad
\Ga_c\times (0,T),\label{eq3.29}\\
&& v_t=0\quad\mbox{on}\quad \Ga_f\times (0,T).\label{eq3.30}
\eeq

\begin{lem}\label{lem3.5}
Let the assumptions in Theorem \ref{th1.1} hold. The following energy inequality holds for $t\in[0,T]$
\begin{equation}\label{eq3.31}
V_1(t)+\int_0^{t}D_1(\tau)d\tau\leqslant V_1(0)+\mid\int^t_0(R_1(\tau),v_t(\tau))d\tau\mid+C\int^t_0L(\tau)d\tau,
\end{equation}
where
\begin{eqnarray}\label{eq3.32}
\int^t_0(R_1(\tau),v_t(\tau))d\tau=&-&\int_0^{t}\int_{\Omega_f}\langle\partial_t(\mathbf{a}:\mathbf{a}^{\mathbf{T}}):Dv,
Dv_{t}\rangle dxd\tau
\nonumber\\
&+&\int_0^{t}\int_{\Omega_f}\langle\partial_t\mathbf{a}q,Dv_{t}\rangle dxd\tau
-\int_0^{t}\int_{\Omega_f}\langle\partial_t\mathbf{a}\partial_tq,Dv\rangle dxd\tau.
\end{eqnarray}
\end{lem}

{\bf Proof} \,\,\, Take $L^2$ inner product with $v_t$ and $w_{tt}$ to \eqref{eq3.25} and \eqref{eq3.27}, respectively. Utilizing the boundary conditions \eqref{eq3.28}-\eqref{eq3.30}, we
attain that
\begin{eqnarray}\label{eq3.33}
&&\frac{1}{2}\frac{d}{dt}\parallel v_t\parallel^2_{L^2}+\int_{\Omega_f}\langle\mathbf{a}:\mathbf{a}^{\mathbf{T}}:Dv_t,Dv_t\rangle dx
+\int_{\Omega_f}\langle\partial_t(\mathbf{a}:\mathbf{a}^{\mathbf{T}}):Dv,Dv_{t}\rangle dx  \\ \nonumber
&&-\int_{\Omega_f}\langle\partial_t(\mathbf{a}q),Dv_{t}\rangle dx +\int_{\Gamma_c}\langle(\omega_t)_{\nu_\mathcal{B}},v_t\rangle d\sigma=0
\end{eqnarray}
and
\begin{eqnarray}\label{eq3.34}
&&\frac{1}{2}\frac{d}{dt}(\parallel w_{tt}\parallel^2_{L^2}+\parallel w_{t}\parallel^2_{L^2}+\int_{\Omega_e}\<l_{Dw_t}D^2W,Dw_t\>dx)
-\frac{1}{2}\int_{\Om_e}D^3W(Dw_t,Dw_t,Dw_t)dx \nonumber \\
&&-\int_{\Gamma_c}\langle(w_t)_{\nu_{\mathcal{B}}},v_t-\gamma(w_t)_{\nu_{\mathcal{B}}}\rangle d\sigma=0.\quad\quad
\end{eqnarray}
Add \eqref{eq3.33} and \eqref{eq3.34} together and integrate in time from $0$ to $t$.
Note that
\begin{equation*}
\int_{\Om_e}D^3W(Dw_t,Dw_t,Dw_t)dx\leqslant CL(t).
\end{equation*}
Due to the ellipticity of $\mathbf{a}(x,t)$,
we get
\begin{eqnarray}\label{eq3.35}
V_1(t)&+&\int_0^{t}D_1(\tau)d\tau\leqslant
V_1(0)+C\int_0^{t}L(\tau)d\tau-\int^t_0\int_{\Omega_f}\langle\partial_t(\mathbf{a}:\mathbf{a}^{\mathbf{T}}):Dv,
Dv_{t}\rangle dxd\tau \nonumber \\
&+&\int^t_0\int_{\Omega_f}\langle\partial_t(\mathbf{a}q),Dv_{t}\rangle dxd\tau
\end{eqnarray}

Taking advantage of \eqref{eq3.26}, we have
\begin{eqnarray}\label{eq3.36}
\int^t_0\int_{\Omega_f}\langle\partial_t(\mathbf{a}q),Dv_{t}dxd\tau
=\int^t_0\int_{\Omega_f}\langle\partial_t\mathbf{a}q,Dv_{t}\rangle dxd\tau
+\int^t_0\int_{\Omega_f}\langle \mathbf{a}\partial_{t}q,Dv_{t}\rangle dxd\tau  \\ \nonumber
=\int^t_0\int_{\Omega_f}\langle\partial_t\mathbf{a}q,Dv_{t}\rangle dxd\tau
-\int^t_0\int_{\Omega_f}\langle\partial_t\mathbf{a}\partial_tq,Dv\rangle dxd\tau
\end{eqnarray}

We submit \eqref{eq3.36} into \eqref{eq3.35} and obtain (\ref{eq3.31}).  \hfill$\Box$\\

According to Lemma \ref{lem3.3}, we find that
\beq\label{eq3.37}
&&\int_{0}^tV_1(\tau)d\tau
\leqslant CV_1(0)+CV_1(t)+C\int_{0}^tL(\tau)d\tau
+C\int_{0}^t\int_{\Ga_c}\mid(w_t)_{\nu_{\mathcal{B}}}\mid^2d\sigma d\tau \nonumber \\
&&+C\int_{0}^t\|Dv_t\|^2_{L^2}d\tau.
\eeq
After a similar procedure with subsection 3.1, we conclude that
\begin{lem}\label{lem3.6}
Under the same hypotheses as in Lemma \ref{eq3.5}, we have for all $t\in[0,T]$ that
\beq\label{eq3.38}
V_1(t)+\int_0^tV_1(\tau)d\tau\leqslant CV_1(0)+C\int_0^tL(\tau)d\tau+C\mid\int^t_0(R_1(\tau),v_t(\tau))d\tau\mid.
\eeq
\end{lem}

\subsection{Third level estimates}
\hskip\parindent
Here, we go further for the third level energy estimates. The third level energy is defined by
\begin{equation*}
V_2(t)=V_2^e(t)+\frac{1}{2}\|v_{tt}\|^2_{L^2},
\end{equation*}
with
\begin{equation*}
V_2^e(t)=\frac{1}{2}(\|w_{ttt}\|^2_{L^2}+\|w_{tt}\|^2_{L^2}+\int_{\Om_e}\langle l_{Dw_{tt}}D^2W,Dw_{tt}\rangle dx)
\end{equation*}
and the dissipative term
\begin{equation*}
D_2(t)=\frac{1}{C}\parallel Dv_{tt}(t)\parallel^2_{L^2}
+\gamma\parallel (w_{tt})_{\nu_{\mathcal{B}}}(t)\parallel^2_{L^2(\Gamma_c)}.
\end{equation*}

Before the derivation of our estimates, we give the equations satisfied by $v_{tt}$ and $w_{tt}.$
\begin{eqnarray}
&&\partial_{t}v_{tt}-\partial_{tt}\div(\mathbf{a}:\mathbf{a}^{\mathbf{T}}:Dv)+\partial_{tt}\div(\mathbf{a}q)=0 \quad \mbox{in}\quad \Omega_f\times (0,T),\label{eq3.39}\\
&&tr(\mathbf{a}_{tt}:Dv)+2tr(\mathbf{a}_{t}:Dv_t)+tr(\mathbf{a}:Dv_{tt})=0  \quad \mbox{in}\quad \Omega_f\times (0,T),\label{eq3.40}\\
&&w_{tttt}-\mathcal{B}_{w}(t)w_{tt}+w_{tt}-r_1(t)=0 \quad \mbox{in}\quad \Omega_e\times (0,T).\label{eq3.41}
\end{eqnarray}

Moreover, the boundary conditions with respect to (\ref{eq3.39})-(\ref{eq3.41}) are the following:
\beq
&&w_{ttt}=v_{tt}-\gamma[(w_{tt})_{\nu_{\mathcal{B}}}+r_{1,\Ga_c}] \quad\mbox{on}\quad \Ga_c\times (0,T),\label{eq3.42}\\
&&(w_{tt})_{\nu_{\mathcal{B}}}=\pl_{tt}(\mathbf{a}:\mathbf{a}^{\mathbf{T}}:Dv)\nu-\pl_{tt}(q\mathbf{a})\nu-r_{1,\Ga_c} \quad\mbox{on}\quad
\Ga_c\times (0,T),\label{eq3.43}\\
&& v_{tt}=0\quad\mbox{on}\quad \Ga_f\times (0,T).\label{eq3.44}
\eeq

\begin{lem}\label{lem3.7}
Suppose that the assumptions in Theorem \ref{th1.1} hold, then the following energy inequality holds for $t\in[0,T]$
\begin{eqnarray}\label{eq3.45}
V_2(t)&+&\int_0^{t}D_2(\tau)d\tau\leqslant CV_2(0)+\mid\int^t_0(R_2(\tau),v_{tt})d\tau\mid+C\int^t_0L(\tau)d\tau \nonumber \\
&+&2\bar{\epsilon}\int_0^t\int_{\Omega_e}\mid w_{ttt}\mid^2dxd\tau,
\end{eqnarray}
where $0<\bar{\epsilon}<1$ sufficiently small and
\begin{eqnarray}\label{eq3.46}
&&\int^t_0(R_2(\tau),v_{tt})d\tau=
2\int_0^{t}\int_{\Omega_f}\langle\partial_t(\mathbf{a}:\mathbf{a}^{\mathbf{T}}):Dv_t,Dv_{tt}\rangle dxd\tau
-\int_0^{t}\int_{\Omega_f}\<\partial_{tt}(\mathbf{a}q),Dv_{tt}\>dxd\tau
\nonumber \\
&&+\int_0^{t}\int_{\Omega_f}\<\partial_{tt}(\mathbf{a}:\mathbf{a}^{\mathbf{T}}):Dv,Dv_{tt}\>dxd\tau.
\end{eqnarray}
\end{lem}

{\bf Proof} \,\,\, Take Euclidean dot product to \eqref{eq3.39} with $v_{tt}$ and integrate over $\Omega_f$. After integrating by parts,
we obtain
\begin{eqnarray}\label{eq3.47}
&&\frac{1}{2}\frac{d}{dt}\parallel v_{tt}\parallel^2_{L^2}
+\int_{\Omega_f}\langle\mathbf{a}:\mathbf{a}^{\mathbf{T}}:Dv_{tt},Dv_{tt}\rangle dx
+\int_{\Gamma_c}\langle(w_{tt})_{\nu_{\mathcal{B}}}+r_{1,\Ga_c},v_{tt}\rangle d\sigma  \nonumber \\
&&+\int_{\Omega_f}\langle\pl_{tt}(\mathbf{a}:\mathbf{a}^{\mathbf{T}}):Dv,Dv_{tt}\rangle dx
+2\int_{\Omega_f}\langle\pl_{t}(\mathbf{a}:\mathbf{a}^{\mathbf{T}}):Dv_t,Dv_{tt}\rangle dx \nonumber \\
&&-\int_{\Omega_f}\langle\pl_{tt}(\mathbf{a}q),Dv_{tt}\rangle dx=0.
\end{eqnarray}
Next, we do the same operation to \eqref{eq3.41} with $w_{ttt}$ as above and also integrate by parts over $\Omega_e$. Thus it follows that
\begin{eqnarray}\label{eq3.48}
&&\frac{1}{2}\frac{d}{dt}(\|w_{ttt}\|^2_{L^2}+\parallel w_{tt}\parallel^2_{L^2}+\int_{\Omega_e}\<l_{Dw_{tt}}D^2W,Dw_{tt}\>dx)
-\int_{\Omega_e}\langle r_1(t),w_{ttt}\rangle dx \nonumber \\
&&-\frac{1}{2}\int_{\Omega_e}D^3W(Dw_{tt},Dw_{tt},Dw_{t})dx
-\int_{\Gamma_c}\langle(w_{tt})_{\nu_{\mathcal{B}}},v_{tt}-\gamma[(w_{tt})_{\nu_{\mathcal{B}}}+r_{1,\Ga_c}]\rangle d\sigma=0.
\end{eqnarray}
Adding \eqref{eq3.47} to \eqref{eq3.48} and integrating in time from $0$ to $t$, it leads to
\begin{eqnarray}\label{eq3.49}
&&V_2(t)+\int_0^t\int_{\Omega_f}\langle\mathbf{a}:\mathbf{a}^{\mathbf{T}}:Dv_{tt},Dv_{tt}\rangle dxd\tau
+\gamma\int_0^t\parallel(w_{tt})_{\nu_{\mathcal{B}}}\parallel^2_{L^2(\Gamma_c)}d\tau  \nonumber \\
&&+\int_0^t\int_{\Gamma_c}\langle r_{1,\Ga_c},v_{tt}\rangle d\sigma d\tau
+\gamma\int_0^t\int_{\Gamma_c}\langle (w_{tt})_{\nu_{\mathcal{B}}},r_{1,\Ga_c}\rangle d\sigma d\tau
+\int^t_0(R_2(\tau),v_{tt})d\tau \nonumber \\
&&=V_2(0)+\int_0^t\int_{\Omega_e}\langle r_1(t),w_{ttt}\rangle dxd\tau
+\frac{1}{2}\int_0^t\int_{\Omega_e}D^3W(Dw_{tt},Dw_{tt},Dw_{t})dxd\tau.
\end{eqnarray}

By Lemma \ref{lem2.5} and the Poincar\'{e} inequality, we have
\beq\label{eq3.50}
\mid\int_0^t\int_{\Gamma_c}\langle r_{1,\Ga_c},v_{tt}\rangle d\sigma d\tau\mid&\leqslant&
C_{\epsilon}\int_0^t\int_{\Gamma_c}\mid r_{1,\Ga_c}\mid^2d\sigma d\tau+\epsilon\int_0^t\int_{\Gamma_c}\mid v_{tt}\mid^2d\sigma d\tau \nonumber \\
&\leqslant&C_{\bar{\epsilon}}\int_0^tL(\tau)d\tau+\bar{\epsilon}\int_0^t\mid Dv_{tt}\mid^2d\tau,
\eeq
and
\be\label{eq3.51}
\gamma\mid\int_0^t\int_{\Gamma_c}\langle(w_{tt})_{\nu_{\mathcal{B}}},r_{1,\Ga_c}\rangle d\sigma d\tau\mid
\leqslant \bar{\epsilon}\int_0^t\int_{\Gamma_c}\mid(w_{tt})_{\nu_{\mathcal{B}}}\mid^2d\sigma d\tau+C_{\bar{\epsilon},\gamma}\int_0^tL(\tau)d\tau,
\ee
where $0<\bar{\epsilon}<1$ is small enough and to be determined.

Similarly, also by using Lemma \ref{lem2.5}, we obtain
\be\label{eq3.52}
\mid\int_0^t\int_{\Omega_e}\langle r_1(t),w_{ttt}\rangle dxd\tau\mid\leqslant \bar{\epsilon}\int_0^t\int_{\Omega_e}\mid w_{ttt}\mid^2dxd\tau
+C_{\bar{\epsilon}}\int_0^tL(\tau)d\tau.
\ee
Note that
\be\label{eq3.53}
\mid\int_0^t\int_{\Omega_e}D^3W(Dw_{tt},Dw_{tt},Dw_{t})dxd\tau\mid\leqslant C\int_0^tL(\tau)d\tau.
\ee
Substitute \eqref{eq3.50}-\eqref{eq3.53} into \eqref{eq3.49} and set $1-C\bar{\epsilon}\geqslant\frac{1}{2}$ and $1-\frac{\bar{\epsilon}}{\gamma}\geqslant\frac{1}{2}$. Thus, via the uniformly
ellipticity, we finally get \eqref{eq3.45}.  \hfill$\Box$\\

From Lemma \ref{lem3.3}, we deduce that
\beq\label{eq3.54}
&&\int_{0}^tV_2(\tau)d\tau
\leqslant CV_2(0)+CV_2(t)+C\int_{0}^tL(\tau)d\tau
+C\int_{0}^t\int_{\Ga_c}\mid(w_{tt})_{\nu_{\mathcal{B}}}\mid^2d\sigma d\tau \nonumber \\
&&+C\int_{0}^t\|Dv_{tt}\|^2_{L^2}d\tau.
\eeq
Multiply \eqref{eq3.54} by $\epsilon'>0$ with $2\bar{\epsilon}<\epsilon'<\min\{\frac{1}{2C},\frac{\gamma}{C}\}$ and add the resulted
inequality to \eqref{eq3.45}. Hence, it turns out that we arrive at the following lemma.
\begin{lem}\label{lem3.8}
Let the hypotheses in Theorem \ref{th1.1} be true. Then we have for $t\in[0,T]$,
\be\label{eq3.55}
V_2(t)+\int_0^{t}V_2(\tau)d\tau\leqslant CV_2(0)+C\mid\int^t_0(R_2(\tau),v_{tt})d\tau\mid+C\int^t_0L(\tau)d\tau.
\ee
\end{lem}

\subsection{Fourth level estimates}
\hskip\parindent
We move on to the Fourth level energy estimates and repeat what we do as above. The fourth level energy is defined by
\begin{equation*}
V_3(t)=V_3^e(t)+\frac{1}{2}\|v_{ttt}\|^2_{L^2} \quad\mbox{and}\quad V_3^e(t)=\frac{1}{2}(\|w_{tttt}\|^2_{L^2}+\|w_{ttt}\|^2_{L^2}+\int_{\Om_e}\langle l_{Dw_{ttt}}D^2W,Dw_{ttt}\rangle dx)
\end{equation*}
as in the previous section. Besides, the dissipative term for the fourth energy is
\begin{equation*}
D_3(t)=\frac{1}{C}\parallel Dv_{ttt}(t)\parallel^2_{L^2}
+\gamma\parallel (w_{ttt})_{\nu_{\mathcal{B}}}(t)\parallel^2_{L^2(\Gamma_c)}.
\end{equation*}

First of all, as before, we differentiate the whole system three times in time and obtain
\begin{eqnarray}
&&\partial_{t}v_{ttt}-\partial_{ttt}\div(\mathbf{a}:\mathbf{a}^{\mathbf{T}}:Dv)+\partial_{ttt}\div(\mathbf{a}q)=0 \quad \mbox{in}\quad \Omega_f\times (0,T),\label{eq3.56}\\
&&tr[\partial_{ttt}(\mathbf{a}:Dv)]=0  \quad \mbox{in}\quad \Omega_f\times (0,T),\label{eq3.57}\\
&&w^{(3)}_{tt}-\mathcal{B}_{w}(t)w_{ttt}+w_{ttt}-r_2(t)=0 \quad \mbox{in}\quad \Omega_e\times (0,T).\label{eq3.58}
\end{eqnarray}

And, the boundary conditions satisfied by the system (\ref{eq3.56})-(\ref{eq3.58}) are as follows:
\beq
&&w_{tttt}=v_{ttt}-\gamma[(w_{ttt})_{\nu_{\mathcal{B}}}+r_{2,\Ga_c}] \quad\mbox{on}\quad \Ga_c\times (0,T),\label{eq3.59}\\
&&(w_{ttt})_{\nu_{\mathcal{B}}}=\pl_{ttt}(\mathbf{a}:\mathbf{a}^{\mathbf{T}}:Dv)\nu-\pl_{ttt}(q\mathbf{a})\nu-r_{2,\Ga_c} \quad\mbox{on}\quad
\Ga_c\times (0,T),\label{eq3.60}\\
&& v_{ttt}=0\quad\mbox{on}\quad \Ga_f\times (0,T).\label{eq3.61}
\eeq

Hence, we are ready to derive the energy estimates for the fourth order energy.

\begin{lem}\label{lem3.9}
Assume that the hypotheses of Theorem \ref{th1.1} hold, then the following energy estimate is true for $t\in[0,T]$
\begin{eqnarray}\label{eq3.62}
V_3(t)&+&\int_0^{t}D_3(\tau)d\tau\leqslant CV_3(0)+\mid\int^t_0(R_3(\tau),v_{ttt})d\tau\mid+C\int^t_0L(\tau)d\tau \nonumber \\
&+&2\epsilon\int_0^tV_3(\tau)d\tau,
\end{eqnarray}
where $0<\tilde{\epsilon}<1$ is sufficiently small and
\begin{eqnarray}\label{eq3.63}
&&\int^t_0(R_3(\tau),v_{ttt})d\tau=\int_0^{t}\int_{\Omega_f}\<\partial_{ttt}(\mathbf{a}:\mathbf{a}^{\mathbf{T}}):Dv,Dv_{ttt}\>dxd\tau
-\int_0^{t}\int_{\Omega_f}\<\partial_{ttt}(\mathbf{a}q),Dv_{ttt}\>dxd\tau
 \nonumber \\
&&+3\int_0^{t}\int_{\Omega_f}\langle\partial_{t}(\mathbf{a}:\mathbf{a}^{\mathbf{T}}):Dv_{tt},Dv_{ttt}\rangle dxd\tau
+3\int_0^{t}\int_{\Omega_f}\langle\partial_{tt}(\mathbf{a}:\mathbf{a}^{\mathbf{T}}):Dv_t,Dv_{ttt}\rangle dxd\tau.
\end{eqnarray}
\end{lem}

{\bf Proof} \,\,\, Take $L^2$ inner product with $v_{ttt}$ and $w^{(4)}$ to \eqref{eq3.56} and \eqref{eq3.58}, respectively.
From \eqref{eq3.59} and \eqref{eq3.60}, we attain that
\beq\label{eq3.64}
&&\frac{1}{2}\frac{d}{dt}\parallel v_{ttt}\parallel^2_{L^2}
+\int_{\Omega_f}\langle\mathbf{a}:\mathbf{a}^{\mathbf{T}}:Dv_{ttt},Dv_{ttt}\rangle dx
+\int_{\Gamma_c}\langle(w_{ttt})_{\nu_{\mathcal{B}}}+r_{2,\Ga_c},v_{ttt}\rangle d\sigma  \nonumber \\
&&+\int_{\Omega_f}\langle\pl_{ttt}(\mathbf{a}:\mathbf{a}^{\mathbf{T}}):Dv,Dv_{ttt}\rangle dx
+3\int_{\Omega_f}\langle\pl_{tt}(\mathbf{a}:\mathbf{a}^{\mathbf{T}}):Dv_t,Dv_{ttt}\rangle dx \nonumber \\
&&+3\int_{\Omega_f}\langle\pl_{t}(\mathbf{a}:\mathbf{a}^{\mathbf{T}}):Dv_{tt},Dv_{ttt}\rangle dx
-\int_{\Omega_f}\langle\pl_{ttt}(\mathbf{a}q),Dv_{ttt}\rangle dx=0
\eeq
and
\begin{eqnarray}\label{eq3.65}
&&\frac{1}{2}\frac{d}{dt}(\|w_{tttt}\|^2_{L^2}+\parallel w_{ttt}\parallel^2_{L^2}+\int_{\Omega_e}\<l_{Dw_{ttt}}D^2W,Dw_{ttt}\>dx) \nonumber \\
&&-\int_{\Omega_e}\langle r_2(t),w_{tttt}\rangle dx
-\frac{1}{2}\int_{\Omega_e}D^3W(Dw_{ttt},Dw_{ttt},Dw_{t})dx \nonumber \\
&&-\int_{\Gamma_c}\langle(w_{ttt})_{\nu_{\mathcal{B}}},v_{ttt}-\gamma[(w_{ttt})_{\nu_{\mathcal{B}}}+r_{2,\Ga_c}]\rangle d\sigma=0.
\end{eqnarray}

Add \eqref{eq3.64} and \eqref{eq3.65} together and we have analogous estimates to \eqref{eq3.50}-\eqref{eq3.53} as well.
Using the similar method with that in Lemma \ref{lem3.7}, we may arrive at \eqref{eq3.62} and conclude the proof.  \hfill$\Box$\\

As a consequence of Lemma \ref{lem3.3}, we have
\beq\label{eq3.66}
&&\int_{0}^tV_3(\tau)d\tau
\leqslant CV_3(0)+CV_3(t)+C\int_{0}^tL(\tau)d\tau
+C\int_{0}^t\int_{\Ga_c}\mid(w_{ttt})_{\nu_{\mathcal{B}}}\mid^2d\sigma d\tau \nonumber \\
&&+C\int_{0}^t\|Dv_{ttt}\|^2_{L^2}d\tau.
\eeq

After the same procedure as that in Subsection 3.3, we acquire the following lemma.
\begin{lem}\label{lem3.10}
Suppose that the hypotheses in Theorem \ref{th1.1} are true. Then for $t\in[0,T]$,
\be\label{eq3.67}
V_3(t)+\int_0^{t}V_3(\tau)d\tau\leqslant CV_3(0)+C\mid\int^t_0(R_3(\tau),v_{ttt})d\tau\mid+C\int^t_0L(\tau)d\tau.
\ee
\end{lem}

\subsection{Superlinear estimates}
\hskip\parindent
Our aim of this subsection is to deal with the perturbation terms in the second, third and fourth level energy estimates. They are
\begin{equation*}
\int^t_0(R_1(\tau),v_t(\tau))d\tau,\quad\int^t_0(R_2(\tau),v_{tt}(\tau))d\tau
\end{equation*}
and
\begin{equation*}
\int^t_0(R_3(\tau),v_{ttt}(\tau))d\tau
\end{equation*}
The concrete presentation of the above three perturbation terms can be found in (\ref{eq3.32}), \eqref{eq3.46} and \eqref{eq3.63}, respectively.
For the estimates of \eqref{eq3.32} and \eqref{eq3.46}, we only list the results. For detail, refer to \cite{IKLT2}.

\begin{lem}$\cite[Lemma \,4.10]{IKLT2}$\label{lem3.11}
We have
\begin{equation*}
\mid(R_1(t),v_t)\mid\leqslant C\parallel v\parallel_{H^1}^{\frac{1}{2}}\parallel v\parallel_{H^2}^{\frac{1}{2}}\parallel v_t\parallel_{H^1}(\parallel v\parallel_{H^2}+\parallel q\parallel_{H^1})+C\parallel v\parallel_{H^1}^{\frac{3}{2}}\parallel v\parallel_{H^2}^{\frac{1}{2}}\parallel q_t\parallel_{H^1},
\end{equation*}
for all $t\in[0,T]$.
\end{lem}

\begin{lem}$\cite[Lemma\, 4.11]{IKLT2}$\label{lem3.12}
For $\epsilon_0\in(0,\frac{1}{C}]$, we have
\begin{eqnarray*}
\mid\int_0^t(R_2(s),v_{tt}(s))ds\mid&\leqslant& \epsilon_0\int_0^t\parallel \nabla v_{tt}\parallel_{L^2}^{2}ds
+C_{\epsilon_0}\int_0^t\parallel v\parallel_{H^1}^{\frac{3}{2}}\parallel v\parallel_{H^3}^{\frac{1}{2}}\parallel q_t\parallel_{H^1}^2ds \\
&+&C_{\epsilon_0}\int_0^t(\parallel v\parallel_{H^3}^{2}+\parallel q\parallel_{H^2}^{2})(\parallel v\parallel_{H^1}^{\frac{5}{2}}\parallel v\parallel_{H^3}^{\frac{3}{2}}+\parallel v_t\parallel_{H^1}^2)ds \\ \nonumber
&+&\epsilon_0\parallel v(t)\parallel_{H^3}^2+\epsilon_0\parallel q_t(t)\parallel_{H^1}^2+\epsilon_0\parallel v_t(t)\parallel_{H^2}^2 \\ \nonumber
&+&C_{\epsilon_0}\parallel v(t)\parallel_{H^1}^{6}\parallel v(t)\parallel_{H^2}^{4}+C_{\epsilon_0}\parallel v(t)\parallel_{H^1}^{2}\parallel v(t)\parallel_{H^2}^{2}\parallel v_t(t)\parallel_{L^2}^2   \\ \nonumber
&+&C\int_0^t(\parallel v\parallel_{H^2}^2+\parallel v_t\parallel_{H^1}^{\frac{1}{2}}\parallel v_t\parallel_{H^2}^{\frac{1}{2}})\parallel q_t\parallel_{H^1}\parallel v_t\parallel_{H^1}ds   \\ \nonumber
&+&C\int_0^t(\parallel v\parallel_{H^2}^3+\parallel v_t\parallel_{H^1}\parallel v\parallel_{H^1}^{\frac{1}{4}}\parallel v\parallel_{H^3}^{\frac{3}{4}})\parallel q_t\parallel_{H^1}\parallel v\parallel_{H^1}^{\frac{3}{4}}\parallel v\parallel_{H^3}^{\frac{1}{4}}ds   \\ \nonumber
&+&C\parallel v(0)\parallel_{H^3}^{6}+C\parallel v_t(0)\parallel_{H^1}^{4}+C\parallel q_t(0)\parallel_{H^1}^2,
\end{eqnarray*}
for all $t\in[0,T]$.
\end{lem}

Now, we turn to \eqref{eq3.63}, even though the computation is quite involved.
\begin{lem}\label{lem3.13}
For $\epsilon_0\in(0,\frac{1}{C}]$, it follows that for all $t\in[0,T]$
\begin{eqnarray*}
&&\int^t_0(R_3(\tau),v_{ttt}(\tau))d\tau\leqslant\epl_0\int^t_0\|Dv_{ttt}\|^2_{L^2}d\tau
+C_{\epl_0}\int^t_0\|v\|^4_{H^3}(\|v\|_{H^1}\cdot\|v\|_{H^3}+\|v_t\|_{H^1} )^2d\tau  \\
&&+C_{\epl_0}\int^t_0(\|v\|^2_{H^3}+\|q\|^2_{H^2})(\|v\|^3_{H^2}+\|v_t\|_{H^1}\|v\|_{H^3}+\|v_{tt}\|_{H^1})^2d\tau  \\
&&+C_{\epl_0}\int^t_0(\|v_t\|^2_{H^2}+\|q_t\|^2_{H^1})(\|v\|^2_{H^2}+\|v_t\|^{\frac{1}{2}}_{H^1}\|v_t\|^{\frac{1}{2}}_{H^2})^2d\tau
\nonumber \\
&&+C_{\epl_0}\int^t_0\|v\|^2_{H^3}(\|v_{tt}\|^2_{H^1}+\|q_{tt}\|^2_{L^2})d\tau
+C_{\epl_0}\int^t_0\|q_{tt}\|^2_{H^1}\|v\|^{\frac{3}{2}}_{H^1}\|v\|^{\frac{1}{2}}_{H^3}d\tau  \\
&&+C\int^t_0(\|v\|^{\frac{7}{2}}_{H^2}\|v\|^{\frac{1}{2}}_{H^3}+\|v\|^{\frac{1}{2}}_{H^2}\|v\|^{\frac{1}{2}}_{H^3}\|v_{tt}\|_{H^1}
+\|v_t\|^2_{H^2})\|q_{tt}\|_{H^1}\|v\|^{\frac{3}{4}}_{H^1}\|v\|^{\frac{1}{4}}_{H^3}d\tau \\
&&+C\int^t_0(\|v\|^{3}_{H^2}+\|v\|^{\frac{1}{2}}_{H^2}\|v\|^{\frac{1}{2}}_{H^3}\|v_{t}\|_{H^1}
+\|v_{tt}\|_{H^1})\|q_{tt}\|_{H^1}\|v_t\|^{\frac{3}{4}}_{H^1}\|v_t\|^{\frac{1}{4}}_{H^3}d\tau \\
&&+C\int^t_0\|v_{tt}\|_{H^1}\|q_{tt}\|_{H^1}(\|v\|^{2}_{H^2}+\|v_t\|^{\frac{3}{4}}_{H^1}\|v_t\|^{\frac{1}{4}}_{H^3})d\tau
 +\epl_0\|q_{tt}\|^2_{H^1}+\epl_0\|v_t\|_{H^3}^2 \nonumber \\
&&+C_{\epl_0}\|v_t\|_{H^1}^6+C\|v\|_{H^3}^4
+C_{\epl_0}\|v_t\|_{H^1}^3\|v\|_{H^1}^6+C_{\epl_0}\|v\|_{H^2}^7\|v\|_{H^1} \nonumber \\
&&+C_{\epl_0}\|v_t\|_{H^1}^2\|v\|_{H^3}^{\frac{3}{2}}\|v\|_{H^1}^{\frac{3}{2}}\|v\|_{H^2}
+C_{\epl_0}\|v_{tt}\|_{H^1}^2\|v\|_{H^1}\|v\|_{H^2}+C\|q_{tt}(0)\|^2_{H^1} \nonumber \\
&&+C\|v_{tt}(0)\|^4_{H^1}+C\|v(0)\|_{H^2}^4+C\|v_t(0)\|_{H^3}^4+C\|v(0)\|_{H^3}^8
\end{eqnarray*}
\end{lem}

{\bf Proof} \,\,\, From \eqref{eq3.63}, we have
\begin{eqnarray}\label{eq3.68}
&&\mid\int^t_0(R_3(\tau),v_{ttt})d\tau\mid\leqslant
\mid\int_0^{t}\int_{\Omega_f}\<\partial_{ttt}(\mathbf{a}:\mathbf{a}^{\mathbf{T}}):Dv,Dv_{ttt}\>dxd\tau\mid \nonumber \\
&&+3\mid\int_0^{t}\int_{\Omega_f}\langle\partial_{t}(\mathbf{a}:\mathbf{a}^{\mathbf{T}}):Dv_{tt},Dv_{ttt}\rangle dxd\tau\mid \nonumber \\
&&+3\mid\int_0^{t}\int_{\Omega_f}\langle\partial_{tt}(\mathbf{a}:\mathbf{a}^{\mathbf{T}}):Dv_t,Dv_{ttt}\rangle dxd\tau\mid \nonumber \\
&&+\mid\int_0^{t}\int_{\Omega_f}\<\partial_{ttt}\mathbf{a}q,Dv_{ttt}\>dxd\tau\mid
+\mid\int_0^{t}\int_{\Omega_f}\<\mathbf{a}q_{ttt},Dv_{ttt}\>dxd\tau\mid  \nonumber \\
&&+3\mid\int_0^{t}\int_{\Omega_f}\<\partial_{tt}\mathbf{a}q_t,Dv_{ttt}\>dxd\tau\mid
+3\mid\int_0^{t}\int_{\Omega_f}\<\partial_{t}\mathbf{a}q_{tt},Dv_{ttt}\>dxd\tau\mid  \nonumber \\
&&=R_{31}+R_{32}+R_{33}+R_{34}+R_{35}+R_{36}+R_{37}
\end{eqnarray}

By H\"{o}lder inequality and Lemma \ref{lem3.1}, we get
\beq\label{eq3.69}
&&\sum_{j=1}^4R_{3j}\leqslant C\int_0^{t}\|Dv_{ttt}\|_{L^2}\|Dv\|_{L^{\infty}}(\|Dv\|_{H^1}^3+\|Dv_t\|_{L^2}\|Dv\|_{L^{\infty}}+ \|Dv_{tt}\|_{L^2})d\tau \nonumber \\
&&+C\int_0^{t}\|Dv_{ttt}\|_{L^2}\|Dv\|^2_{L^{\infty}}(\|Dv\|_{L^2}\|Dv\|_{L^{\infty}}+\|Dv_t\|_{L^2})d\tau \nonumber \\
&&+C\int_0^{t}\|Dv_{ttt}\|_{L^2}\|Dv_{tt}\|_{L^2}\|Dv\|_{L^{\infty}}d\tau
+C\int_0^{t}\|Dv_{ttt}\|_{L^2}\|Dv_{t}\|_{L^2}\|Dv\|^2_{L^{\infty}}d\tau \nonumber \\
&&+C\int_0^{t}\|Dv_{ttt}\|_{L^2}\|Dv_{t}\|_{L^6}(\|v\|_{H^2}^2+\|Dv_t\|_{L^3})d\tau \nonumber \\
&&+C\int_0^{t}\|q\|_{L^{\infty}}\|Dv_{ttt}\|_{L^2}(\|v\|^3_{H^2}+\|v_t\|_{H^1}\|v\|_{H^3}+\|v_{tt}\|_{H^1})d\tau \nonumber \\
&&\leqslant C\int_0^{t}\|Dv_{ttt}\|_{L^2}\|v\|_{H^3}(\|v\|_{H^2}^3+\|v_t\|_{H^1}\|v\|_{H^3}+ \|v_{tt}\|_{H^1})d\tau \nonumber \\
&&+C\int_0^{t}\|Dv_{ttt}\|_{L^2}\|v\|^2_{H^3}(\|v\|_{H^1}\|v\|_{H^3}+\|v_t\|_{H^1})d\tau \nonumber \\
&&+C\int_0^{t}\|Dv_{ttt}\|_{L^2}\|v_{tt}\|_{H^1}\|v\|_{H^3}d\tau
+C\int_0^{t}\|Dv_{ttt}\|_{L^2}\|v_{t}\|_{H^1}\|v\|^2_{H^3}d\tau \nonumber \\
&&+C\int_0^{t}\|Dv_{ttt}\|_{L^2}\|v_{t}\|_{H^2}(\|v\|_{H^2}^2+\|v_t\|^{\frac{1}{2}}_{H^1}\|v_t\|^{\frac{1}{2}}_{H^2})d\tau \nonumber \\
&&+C\int_0^{t}\|q\|_{H^2}\|Dv_{ttt}\|_{L^2}(\|v\|^3_{H^2}+\|v_t\|_{H^1}\|v\|_{H^3}+\|v_{tt}\|_{H^1})d\tau
\eeq
and
\beq\label{eq3.70}
R_{35}+R_{36}&\leqslant& C\int_0^{t}\|Dv_{ttt}\|_{L^2}\|q_t\|_{H^1}(\|v\|^2_{H^2}+\|v_t\|_{H^1}^{\frac{1}{2}}\|v_t\|_{H^2}^{\frac{1}{2}})d\tau \nonumber \\
&+& C\int_0^{t}\|Dv_{ttt}\|_{L^2}\|q_{tt}\|_{L^2}\|v\|_{H^3}d\tau,
\eeq
where the Sobolev and interpolation inequalities are employed. Now we begin to treat $R_{37}.$ By the differentiated divergence-free
condition \eqref{eq3.57}, we deduce that
\begin{eqnarray*}
R_{37}&=&-\int_0^{t}\int_{\Om_f}\langle q_{ttt}\mathbf{a}_{ttt},Dv\rangle dxd\tau
-3\int_0^{t}\int_{\Om_f}\<q_{ttt}\mathbf{a}_{tt},Dv_t\>dxd\tau
-3\int_0^{t}\int_{\Om_f}\<q_{ttt}\mathbf{a}_{t},Dv_{tt}\>dxd\tau \\
&=&-\int_{\Om_f}\langle q_{tt}\mathbf{a}_{ttt},Dv\rangle dx\mid^t_0-3\int_{\Om_f}\<q_{tt}\mathbf{a}_{tt},Dv_t\>dx\mid^t_0
-3\int_{\Om_f}\<q_{tt}\mathbf{a}_{t},Dv_{tt}\>dx\mid^t_0 \\
&+&\int_0^{t}\int_{\Om_f}\<q_{tt}\pl^4_t\mathbf{a},Dv\>dxd\tau
+4\int_0^{t}\int_{\Om_f}\<q_{tt}\mathbf{a}_{ttt},Dv_t\>dxd\tau \\
&+&6\int_0^{t}\int_{\Om_f}\<q_{tt}\mathbf{a}_{tt},Dv_{tt}\>dxd\tau
+3\int_0^{t}\int_{\Om_f}\<q_{tt}\mathbf{a}_{t},Dv_{ttt}\>dxd\tau.
\end{eqnarray*}

Applying Lemma \ref{lem2.1} and Corollary \ref{cor2.1} along with H\"{o}lder's, Sobolev and interpolation inequalities, we have
\begin{eqnarray*}
&&R_{37}\leqslant \|\mathbf{a}_{ttt}(t)\|_{L^2}\|Dv(t)\|_{L^3}\|q_{tt}(t)\|_{L^6}+3\|\mathbf{a}_{tt}(t)\|_{L^2}\|Dv_{t}(t)\|_{L^3}\|q_{tt}(t)\|_{L^6} \\
&&+3\|\mathbf{a}_{t}(t)\|_{L^3}\|Dv_{tt}(t)\|_{L^2}\|q_{tt}(t)\|_{L^6}+\|\mathbf{a}_{ttt}(0)\|_{L^2}\|Dv(0)\|_{L^3}\|q_{tt}(0)\|_{L^6} \\
&&+3\|\mathbf{a}_{tt}(0)\|_{L^2}\|Dv_{t}(0)\|_{L^3}\|q_{tt}(0)\|_{L^6}
+3\|\mathbf{a}_{t}(0)\|_{L^3}\|Dv_{tt}(0)\|_{L^2}\|q_{tt}(0)\|_{L^6}  \\
&&+\int_0^{t}\|\pl^4_t\mathbf{a}\|_{L^2}\|Dv\|_{L^3}\|q_{tt}\|_{L^6}d\tau
+4\int_0^{t}\|\mathbf{a}_{ttt}\|_{L^2}\|Dv_{t}\|_{L^3}\|q_{tt}\|_{L^6}d\tau \\
&&+6\int_0^{t}\|\mathbf{a}_{tt}(0)\|_{L^3}\|Dv_{tt}\|_{L^2}\|q_{tt}\|_{L^6}d\tau
+3\int_0^{t}\|\mathbf{a}_{t}\|_{L^3}\|Dv_{ttt}\|_{L^2}\|q_{tt}\|_{L^6}d\tau.
\end{eqnarray*}

The sum of the first three terms on the right hand side of the above estimate is bounded by
\begin{eqnarray*}
&&C(\|Dv\|_{H^1}^3+\|Dv_t\|_{L^2}\|Dv\|_{L^{\infty}}
+\|Dv_{tt}\|_{L^2})\|q_{tt}\|_{H^1}\|v\|_{H^1}^{\frac{1}{2}}\|v\|_{H^2}^{\frac{1}{2}} \\
&&+C(\|Dv\|_{L^2}\|Dv\|_{L^{\infty}}+\|D v_t\|_{L^2})\|q_{tt}\|_{H^1}\|v_t\|_{H^1}^{\frac{3}{4}}\|v_t\|_{H^3}^{\frac{1}{4}} \\
&&+C\|Dv\|_{L^3}\|v_{tt}\|_{H^1}\|q_{tt}\|_{H^1}\leqslant \epl_0\|q_{tt}\|^2_{H^1}+\epl_0\|v_t\|_{H^3}^2+C_{\epl_0}\|v_t\|_{H^1}^6 \\
&&+C\|v\|_{H^3}^4+C_{\epl_0}\|v_t\|_{H^1}^3\|v\|_{H^1}^6+C_{\epl_0}\|v\|_{H^2}^7\|v\|_{H^1} \\
&&+C_{\epl_0}\|v_t\|_{H^1}^2\|v\|_{H^3}^{\frac{3}{2}}\|v\|_{H^1}^{\frac{3}{2}}\|v\|_{H^2}
+C_{\epl_0}\|v_{tt}\|_{H^1}^2\|v\|_{H^1}\|v\|_{H^2}
\end{eqnarray*}
with the help of Lemma \ref{lem3.1}. Thus, thanks to the Agmon's inequality $\|\phi\|_{L^{\infty}}\leqslant C\|\phi\|_{H^1}^{\frac{1}{2}}\|\phi\|_{H^2}^{\frac{1}{2}}$ in particular, we obtain
\beq\label{eq3.71}
&&R_{37}\leqslant \epl_0\|q_{tt}\|^2_{H^1}+\epl_0\|v_t\|_{H^3}^2+C_{\epl_0}\|v_t\|_{H^1}^6+C\|v\|_{H^3}^4 \nonumber \\
&&+C_{\epl_0}\|v_t\|_{H^1}^3\|v\|_{H^1}^6+C_{\epl_0}\|v\|_{H^2}^7\|v\|_{H^1}
+C_{\epl_0}\|v_t\|_{H^1}^2\|v\|_{H^3}^{\frac{3}{2}}\|v\|_{H^1}^{\frac{3}{2}}\|v\|_{H^2} \nonumber \\
&&+C_{\epl_0}\|v_{tt}\|_{H^1}^2\|v\|_{H^1}\|v\|_{H^2}+C\|q_{tt}(0)\|^2_{H^1}+C\|v_{tt}(0)\|^4_{H^1} \nonumber \\
&&+C\|v(0)\|_{H^2}^4+C\|v_t(0)\|_{H^3}^4+C\|v(0)\|_{H^3}^8+C\int_0^{t}\|q_{tt}\|_{H^1}\|Dv_{ttt}\|_{L^2}
\|v\|_{H^1}^{\frac{3}{4}}\|v\|_{H^3}^{\frac{1}{4}}d\tau \nonumber \\
&&+C\int^t_0(\|v\|^{3}_{H^2}+\|v\|^{\frac{1}{2}}_{H^2}\|v\|^{\frac{1}{2}}_{H^3}\|v_{t}\|_{H^1}
+\|v_{tt}\|_{H^1})\|q_{tt}\|_{H^1}\|v_t\|^{\frac{3}{4}}_{H^1}\|v_t\|^{\frac{1}{4}}_{H^3}d\tau \nonumber \\
&&+C\int^t_0(\|v\|^{\frac{7}{2}}_{H^2}\|v\|^{\frac{1}{2}}_{H^3}+\|v\|^{\frac{1}{2}}_{H^2}\|v\|^{\frac{1}{2}}_{H^3}\|v_{tt}\|_{H^1}
+\|v_t\|^2_{H^2} \nonumber \\
&&+\|Dv_{ttt}\|_{L^2})\|q_{tt}\|_{H^1}\|v\|^{\frac{3}{4}}_{H^1}\|v\|^{\frac{1}{4}}_{H^3}d\tau.
\eeq

Therefore, combining \eqref{eq3.69}-\eqref{eq3.71}, we conclude the proof of this lemma.  \hfill$\Box$\\

\section{Energy decay and global existence of the system}
\setcounter{equation}{0}
\hskip\parindent
We aim at the global existence of solutions and the energy decay estimates in this section.

Let the total energy of the whole system
\begin{equation}\label{total}
X(t)=\sum_{i=0}^3V_i(t)+\epsilon_1(\|Dv\|^2_{L^2}+\|Dv_t\|^2_{L^2}+\|Dv_{tt}\|^2_{L^2})
\end{equation}
and its equivalent version
\begin{equation*}
\mathcal{X}(t)=\frac{1}{2}\sum_{i=0}^3\|v^{(j)}\|^2_{L^2}+\mathcal{E}^e(t)
+\epsilon_1(\|Dv\|^2_{L^2}+\|Dv_t\|^2_{L^2}+\|Dv_{tt}\|^2_{L^2})
\end{equation*}
where $\epsilon_1>0$  is given sufficiently small and to be determined later.

We make some preparations for the proof of Theorem \ref{th1.1}.

We have
\beq\label{eq4.1}
\parallel Dv(t)\parallel^2_{L^2}&&=\parallel Dv(0)\parallel^2_{L^2}+\int^t_0\frac{d}{d\tau}\parallel D v(\tau)\parallel^2_{L^2}d\tau\nonumber\\
&&= \parallel Dv(0)\parallel^2_{L^2}+2\int^t_0\parallel Dv\parallel_{L^2}\parallel Dv_t\parallel_{L^2}d\tau\nonumber\\
&&\leq\parallel Dv(0)\parallel^2_{L^2}+C\int^t_0(D_0(\tau)+D_1(\tau))d\tau,
\eeq
Similarly, we obtain
\begin{equation}\label{eq4.2}
\parallel Dv_t(t)\parallel^2_{L^2}\leqslant\parallel Dv_t(0)\parallel^2_{L^2}+C\int^t_0(D_1(\tau)+D_2(\tau))d\tau
\end{equation}
and
\begin{equation}\label{eq4.3}
\parallel Dv_{tt}(t)\parallel^2_{L^2}\leqslant\parallel Dv_{tt}(0)\parallel^2_{L^2}+C\int^t_0(D_2(\tau)+D_3(\tau))d\tau.
\end{equation}

In addition, it follows from Lemmas \ref{lem3.1} and \ref{lem3.4} that
\begin{equation}\label{eq4.4}
V_0(t)+\int^t_0V_0(\tau)d\tau+\int^t_0D_0(\tau)d\tau\leqslant CV_0(0)+C\int_0^tL(\tau)d\tau.
\end{equation}
Combining Lemmas \ref{lem3.5}, \ref{lem3.6} and \ref{lem3.11}, we have
\beq\label{eq4.5}
&&V_1(t)+\int^t_0V_1(\tau)d\tau+\int^t_0D_1(\tau)d\tau\leqslant CV_1(0)+C\int_0^tL(\tau)d\tau \nonumber \\
&&+C\int^t_0\mathbf{P}_1(\parallel v\parallel_{H^2},\parallel q\parallel_{H^1},\parallel v_t\parallel_{H^1},\parallel q_t\parallel_{H^1})d\tau.
\eeq
From Lemmas \ref{lem3.7}, \ref{lem3.8} and \ref{lem3.12},
\begin{eqnarray}\label{eq4.6}
V_2(t)&+&\int^t_0V_2(\tau)d\tau+\int^t_0D_2(\tau)d\tau\leqslant CV_2(0)+C\int^t_0L(\tau)d\tau+\epsilon_0\parallel v(t)\parallel_{H^3}^2 \nonumber\\
&+&\epsilon_0\parallel q_t(t)\parallel_{H^1}^2
+\epsilon_0\parallel v_t(t)\parallel_{H^2}^2
+\epsilon_0\int_0^t\parallel Dv_{tt}\parallel_{L^2}^{2}d\tau
+\mathbf{P}_2(\parallel v\parallel_{H^2},\parallel v_t\parallel_{L^2})  \nonumber\\
&+&\int^t_0\mathbf{P}_3(\parallel v\parallel_{H^3},\parallel q\parallel_{H^2},\parallel v_t\parallel_{H^2},\parallel q_t\parallel_{H^1})d\tau   \nonumber\\
&+&\mathbf{P}_4(\parallel v(0)\parallel_{H^3},\parallel v_t(0)\parallel_{H^1},\parallel q_t(0)\parallel_{H^1}).
\end{eqnarray}
Moreover, according to Lemma \ref{lem3.9}, \ref{lem3.10} and \ref{lem3.13}, we attain
\beq\label{eq4.7}
&&V_3(t)+\int^t_0V_3(\tau)d\tau+\int^t_0D_3(\tau)d\tau\leqslant CV_3(0)+C\int^t_0L(\tau)d\tau+\epl_0\|q_{tt}\|^2_{H^1} \nonumber \\
&&+\epl_0\|v_t\|_{H^3}^2+\epl_0\int^t_0\|Dv_{ttt}\|^2_{L^2}d\tau
+\mathbf{P}_5(\|v\|_{H^3},\|v_t\|_{H^1},\|v_{tt}\|_{H^1})  \nonumber\\
&&+\int^t_0\mathbf{P}_6(\|v\|_{H^4},\|q\|_{H^3},\|v_t\|_{H^3},\|q_t\|_{H^2},\|v_{tt}\|_{H^2},\|q_{tt}\|_{H^1})d\tau   \nonumber\\
&&+\mathbf{P}_7(\|v(0)\|_{H^3},\|v_t(0)\|_{H^3},\|q_{tt}(0)\|_{H^1},\|v_{tt}(0)\|_{H^1}),
\eeq
where the symbols $\mathbf{P}_i,1\leqslant i\leqslant 7$ denote the superlinear polynomials of
their arguments, which are allowed to depend on $\epsilon_0$ from Lemmas \ref{lem3.12} and \ref{lem3.13}. Now multiply \eqref{eq4.1}-
\eqref{eq4.3} by
sufficiently small $\epsilon_1$, sum up \eqref{eq4.4}-\eqref{eq4.7} and then add them together to obtain
\begin{eqnarray}\label{eq4.8}
&&X(t)+\int^t_0X(\tau)d\tau\leqslant CX(0)+\epl_0\|v(t)\|_{H^3}^2
+\epl_0\|q_t(t)\|_{H^1}^2
+\epl_0\|v_t(t)\|_{H^3}^2  \nonumber\\
&&+\epl_0\|q_{tt}\|^2_{H^1}+C\int^t_0L(\tau)d\tau
+\bar{\mathbf{P}}_1(\parallel v\parallel_{H^3},\parallel v_t\parallel_{H^1},\|v_{tt}\|_{H^1}) \nonumber\\
&&+\int^t_0\bar{\mathbf{P}}_2(\|v\|_{H^4},\|q\|_{H^3},\|v_t\|_{H^3},\|q_t\|_{H^2},\|v_{tt}\|_{H^2},\|q_{tt}\|_{H^1})d\tau  \nonumber \\
&&+\bar{\mathbf{P}}_3(\|v(0)\|_{H^3},\|v_t(0)\|_{H^3},\|q_{tt}(0)\|_{H^1},\|v_{tt}(0)\|_{H^1},\parallel q_t(0)\parallel_{H^1}).
\end{eqnarray}

Because of \eqref{eq2.25} in Theorem \ref{th2.1}, we find that
\beq\label{eq4.9}
&&\mathcal{X}(t)+\int^t_0\mathcal{X}(\tau)d\tau\leqslant C\mathcal{X}(0)+C\epl_0\|v(t)\|_{H^3}^2
+C\epl_0\|q_t(t)\|_{H^1}^2
+C\epl_0\|v_t(t)\|_{H^3}^2  \nonumber\\
&&+C\epl_0\|q_{tt}\|^2_{H^1}+C\int^t_0L(\tau)d\tau
+\bar{\mathbf{P}}_1(\parallel v\parallel_{H^3},\parallel v_t\parallel_{H^1},\|v_{tt}\|_{H^1}) \nonumber\\
&&+\int^t_0\bar{\mathbf{P}}_2(\|v\|_{H^4},\|q\|_{H^3},\|v_t\|_{H^3},\|q_t\|_{H^2},\|v_{tt}\|_{H^2},\|q_{tt}\|_{H^1})d\tau  \nonumber \\
&&+\bar{\mathbf{P}}_3(\|v(0)\|_{H^3},\|v_t(0)\|_{H^3},\|q_{tt}(0)\|_{H^1},\|v_{tt}(0)\|_{H^1},\parallel q_t(0)\parallel_{H^1}) \nonumber\\
&&+CL(t)+C\sum_{j=0}^2\|w^{(j)}_{\nu_{\mathcal{B}}}\|^2_{H^{\frac{5}{2}-j}(\Ga_c)}+
C\sum_{j=0}^2\int^t_0\|w^{(j)}_{\nu_{\mathcal{B}}}\|^2_{H^{\frac{5}{2}-j}(\Ga_c)}d\tau.
\eeq

From \eqref{eq2.4}, Lemma \ref{lem2.4} and the properties of Sobolev space we list, we have
\begin{equation}\label{eq4.10}
\parallel v\parallel_{H^3}^2+\parallel q\parallel_{H^2}^2\leqslant C\mathcal{X}(t).
\end{equation}
Thanks to \eqref{eq2.5} and \eqref{eq4.10}, similarly we obtain
\begin{equation}\label{eq4.11}
\parallel v_t\parallel_{H^3}^2+\parallel q_t\parallel_{H^2}^2\leqslant C\mathcal{X}(t)+C\mathcal{X}^2(t).
\end{equation}

For \eqref{eq2.4} in case of $s=2,$ by \eqref{eq4.11} and Lemma \ref{lem2.4},
\be\label{eq4.12}
\parallel v\parallel_{H^4}^2+\parallel q\parallel_{H^3}^2\leqslant C\mathcal{X}(t)+C\mathcal{X}^2(t).
\ee

From $(2)$ in Lemma \ref{lem2.2}, \eqref{eq4.10} and \eqref{eq4.12}, via a similar way, we deduce that
\begin{equation}\label{eq4.13}
\parallel v_{tt}\parallel_{H^2}^2+\parallel q_{tt}\parallel_{H^1}^2\leqslant C\mathcal{X}(t)+C\tilde{\mathbf{P}}(\mathcal{X}(t)),
\end{equation}
where $\tilde{\mathbf{P}}$ is a polynomial with the degree of each term of it is at least $2.$

Submitting \eqref{eq4.10}-\eqref{eq4.13} and the boundary condition \eqref{eq3.11} into \eqref{eq4.9} and setting
$\epl_0$ small enough and $\gamma\geqslant 2C,$ where the constant $C$ depends on $\mathcal{X}(0),$ by Lemma \ref{lem2.5}, it follows from \eqref{eq4.9}
that
\begin{equation}\label{eq4.14}
\mathcal{X}(t)+\int^t_0\mathcal{X}(\tau)d\tau\leqslant C\mathcal{X}(0)+\mathbf{P}(\mathcal{X}(t))+\int^t_0\mathbf{P}(\mathcal{X}(\tau))d\tau+\mathbf{P}(\mathcal{X}(0)),
\end{equation}
where $\mathbf{P}$ is a superlinear polynomial as well. We rewrite \eqref{eq4.14} as
\begin{equation}\label{eq4.15}
\mathcal{X}(t)+\int^t_{0}\mathcal{X}(\tau)d\tau\leqslant C_0\sum_{j=1}^{m}\int^t_{0}\mathcal{X}(\tau)^{\alpha_j}d\tau+C_0\sum_{k=1}^{n}\mathcal{X}(t)^{\beta_k}
+C_0\sum_{k=1}^{n}\mathcal{X}(0)^{\beta_k}+C_0\mathcal{X}(0),
\end{equation}
where $C_0\geqslant 1$, $\alpha_1,...,\alpha_m>1$ and $\beta_1,...,\beta_n>1$.

Following \cite[Lemma\, 5.1]{IKLT2}, we have
\begin{lem}\label{lem4.1}
Suppose that $\mathcal{X}:[0,\infty)\rightarrow[0,\infty)$ is continuous for all $t$ such that $\mathcal{X}(t)$ is finite and assume that it satisfies
\begin{equation*}
\mathcal{X}(t)+\int^t_{\tau}\mathcal{X}(s)ds\leqslant C_0\sum_{j=1}^{m}\int^t_{\tau}\mathcal{X}(s)^{\alpha_j}ds+C_0\sum_{k=1}^{n}\mathcal{X}(t)^{\beta_k}
+C_0\sum_{k=1}^{n}\mathcal{X}(\tau)^{\beta_k}+C_0\mathcal{X}(\tau),
\end{equation*}
where $\alpha_1,...,\alpha_m>1$ and $\beta_1,...,\beta_n>1.$ Also, assume that $\mathcal{X}(0)\leqslant\epsilon$. If $\epsilon\leqslant\frac{1}{C}$, where the constant $C$ depends on $C_0, m, \alpha_1,...,\alpha_m, \beta_1,...,\beta_n$, we have
$\mathcal{X}(t)\leqslant C\epsilon e^{-\frac{t}{C}}$.
\end{lem}

{\bf Proof of Theorem \ref{th1.1}}\,\,\, Utilizing Lemma \ref{lem4.1} and following the proof of \cite[Theorem 2.1]{IKLT2} and \cite[Theorem 1.1]{Zhang-Yao 2009}, the proof of Theorem \ref{th2.1} is finished.\hfill$\Box$\\

\end{document}